\documentclass{article}

\usepackage{arxiv}

\usepackage[utf8]{inputenc} % allow utf-8 input
\usepackage[T1]{fontenc}    % use 8-bit T1 fonts
\usepackage{hyperref}       % hyperlinks
\usepackage{url}            % simple URL typesetting
\usepackage{booktabs}       % professional-quality tables
\usepackage{amsfonts}       % blackboard math symbols
\usepackage{nicefrac}       % compact symbols for 1/2, etc.
\usepackage{microtype}      % microtypography
\usepackage{lipsum}
\usepackage{graphicx}
\usepackage{color}
\usepackage{mathtools}
\usepackage{cleveref}
\graphicspath{ {./images/} }
\usepackage{subcaption}
\usepackage{verbatim} % command "comment"

\newtheorem{remark}{Remark}[section]

\title{Information-geometry of physics-informed statistical manifolds and its use in data assimilation}

\author{
  Francesca Boso \\
  Department of Energy Resources Engineering\\
  Stanford University\\
  Stanford, CA 94305 \\
  \texttt{fboso@stanford.edu} \\
  \AND
  Daniel M. Tartakovsky \\
  Department of Energy Resources Engineering\\
  Stanford University\\
  Stanford, CA 94305 \\
  \texttt{tartakovsky@stanford.edu} \\
}

\begin{document}
\maketitle
\begin{abstract}
The data-aware method of distributions (DA-MD) is a low-dimension data assimilation procedure to forecast the behavior of dynamical systems described by differential equations. It combines sequential Bayesian update with the MD, such that the former utilizes available observations while the latter propagates the (joint) probability distribution of the uncertain system state(s). The core of DA-MD is the minimization of a distance between an observation and a prediction in distributional terms, with prior and posterior distributions constrained on a statistical manifold  defined by the MD. We leverage the information-geometric properties of the statistical manifold to reduce predictive uncertainty via data assimilation. Specifically, we exploit the information geometric structures induced by two discrepancy metrics, the Kullback-Leibler divergence and the Wasserstein distance, which explicitly yield natural gradient descent. To further accelerate optimization, we build a deep neural network as a surrogate model for the MD that enables automatic differentiation. The manifold's geometry is quantified without sampling, yielding an accurate approximation of the gradient descent direction. Our numerical experiments demonstrate that accounting for the information-geometry of the manifold significantly reduces the computational cost of data assimilation by facilitating the calculation of gradients and by reducing the number of required iterations. Both storage needs and computational cost depend on the dimensionality of a statistical manifold, which is typically small by MD construction. When convergence is achieved, the Kullback-Leibler and $L_2$ Wasserstein metrics have similar performances, with the former being more sensitive to poor choices of the prior.
\end{abstract}

% keywords can be removed
\keywords{Method of Distributions \and Data assimilation \and Uncertainty Reduction \and Machine Learning}

\section{Introduction}

Mathematical models used to represent ``reality'' are invariably faulty due to a number of mutually reinforcing reasons such as lack of detailed knowledge  of the relevant laws of nature, scarcity (in quality and/or quantity) of observations, and inherent spatiotemporal variability of the coefficients used in their parameterizations. Consequently, model predictions must be accompanied by a quantifiable measure of predictive uncertainty (e.g., error bars or confidence intervals); when available, observations should be used to reduce this uncertainty. The probabilistic framework provides a natural means to achieve both goals. 
For example, a random forcing in Langevin (stochastic ordinary-differential) equations~\cite{risken1996fokker} or fluctuating Navier-Stokes (stochastic partial-differential) equations~\cite{landau1980statistical} implicitly account for sub-scale variability and processes that are otherwise absent in the underlying model. %For example, they are routinely used to describe the motion of particles in fluids. A probabilistic representation of such forcings - either generated within the system or externally imposed - and their effect on the model solution allows one to formulate quantitative predictions of system behavior. 

Solutions of such stochastic models, and of models with random coefficients, are given in terms of the (joint) probability density function (PDF) or cumulative distribution function (CDF) of the system state(s). They can be computed, with various degrees of accuracy and ranges of applicability, by employing, e.g., Monte Carlo simulations (MCS), polynomial chaos expansions (PCE) and the method of distributions (MD) \cite{tartakovsky2016method}. MCS are robust, straightforward and trivially parallelizable; yet, they carry (often prohibitively) high computational cost. PCE rely on a finite-dimensional expansion of the solution of a stochastic model; their accuracy and computational efficiency decrease as the correlation length of the random inputs decreases (the so-called curse of dimensionality), making them ill-suited to problems with white noise. The MD yields a (generally approximate) partial differential equation (PDE) for the PDF or CDF of a system state (henceforth referred to as a PDF/CDF equation). The MD can handle inputs with both long and short correlations, although the correlation length might affect the robustness of the underlying closure approximations when the latter are needed. For Langevin systems driven by white noise, the MD yields a Fokker-Planck equation \cite{risken1996fokker} for a system state's PDF. For colored (correlated) noise, PDF/CDF equations become approximate~\cite{wang-2013-probability}, although their computational footprint typically does not change. If a Langevin system is characterized by $N_\text{st}$ system states, then  PDF/CDF equations are defined in an augmented $N_\text{st}$-dimensional space. Their MD-based derivation requires a closure approximation~\cite[and references therein]{tartakovsky2016method} such as the semi-local closure~\cite{barajas2016probabilistic, maltba-2018-nonlocal} used in our analysis because of its accuracy and manageable computational cost. 

The temporal evolution of the PDF of a system state predicted with, e.g., the MD provides a measure of the model's predictive uncertainty in the absence of observations of the system state. In the lingo of Bayesian statistics, this PDF serves as a prior that can be improved (converted into the posterior PDF) via Bayesian update as data become available. %Overestimating the uncertainty affecting the random forcings leads to an inaccurate representation of the uncertainty affecting the state, regardless of the chosen uncertainty propagation technique. Available observations can thus be used to reduce the uncertainty on the state predictions. 
When used in combination with ensemble methods like MCS, standard strategies for Bayesian data assimilation, e.g., Markov chain Monte Carlo (MCMC) and its variants, are often prohibitively expensive \cite{wikle2007bayesian}. The computational expedience is the primary reason for the widespread use of various flavors of Kalman filter, which perform well when the system state's PDF is (nearly) Gaussian and models are linear, but are known to fail otherwise. 
%Alongside classic data assimilation techniques that represent approximations of Bayesian update whenever the latter becomes , e.g. Markov-Chain Monte Carlo (MCMC) and various flavours of filtering, we have recently proposed a novel 
Data-aware MD (DA-MD) \cite{boso-2020-learning} alleviates this  computational bottleneck, rendering Bayesian update feasible even on a laptop. DA-MD employs the MD to propagate the system state PDF (forecast step) and sequential Bayesian update at measurement locations to assimilate data (analysis step). It offers two major benefits. First, the MD replaces repeated model runs, characteristic of both MCMC \cite{brooks2011handbook} and ensemble and particle filters \cite{evensen2009data,arulampalam2002tutorial}, with the solution of a single deterministic equation for the evolving PDF. Second, it dramatically reduces the dimensionality of the PDFs involved in the Bayesian update at each assimilation step because it relies on a single-point PDF rather than a multi-point PDF whose dimensionality is determined by the discretized state being updated. %the number of data points. 
DA-MD takes advantage of the MD's ability to handle nonlinear models and non-Gaussian distributions \cite{boso-2014-cumulative,boso-2020-data}. 

DA-MD recasts data assimilation as a minimization problem, whose loss function represents the discrepancy between observed and predicted posterior distributions. The observed posterior PDF is obtained by direct application of Bayes' rule at the  measurement point, combining the data model and a prior PDF computed via the MD. The predicted PDF is assumed to obey the PDF equation, which acts as a PDE constraint for the loss function. The parameters appearing in the MD are the target of minimization and introduce a suitable parameterization for the space of probabilities (a statistical manifold) with quantifiable geometric properties. The computational effort of DA-MD is thus determined by the efficiency in the solution of a minimization problem on a manifold. This aspect of DA-MD is the central focus of our analysis, in which we exploit information-geometric theory to reformulate the optimization problem by relying on the geometric properties of the MD-defined manifold.

We utilize results from the optimal transport theory and machine learning. 
%This links our work to machine learning and its numerous applications, and to optimal transport theory. %From a different perspective, the method can be seen as data-driven learning of a probabilistic propagation model. physics-informed machine learning or data-driven modeling. We minimize but not in a random space, it has a specific structure that we can take advantage of for both gradient descent and natural gradient descent.
%
Specifically, we employ both the Kullback-Leibler (KL) divergence and the $L_2$ Wasserstein distance to measure the discrepancy between predicted and observed posterior distributions at each assimilation point. The former underpins much of information theory \cite{kullback1997information} and variational inference \cite{blei2017variational}\footnote{Unlike traditional variational inference, our approach utilizes univariate (single-point) distributions that are characterized by a specific, physics-driven parameterization enabled by the MD.}, while the latter has its origins in optimal transport and is now increasingly popular in the wider machine learning community \cite{peyre2019computational}.  We employ gradient descent (GD) and natural gradient descent (NGD) for optimization \cite{chen2018wasserstein}, with preconditioning matrices expressing the geometry induced on the statistical manifold by the choice of the discrepancy. These formulations are explicit for univariate distributions; thus, they ideally suit our data assimilation procedure. 

Finally, we construct a surrogate model for the solution of the PDF/CDF equation to accelerate sequential minimization of loss functions, taking advantage of the relatively small dimensionality of the statistical manifold. We identify a special architecture of a deep neural network (DNN) that enables the calculation of the terms involved in NGD for both discrepancy choices. The use of DNNs obviates the need to resort to sampling when assessing the manifold's geometry, a strategy that has a debatable success~\cite{kunstner2019limitations}. 

The paper is organized as follows. In \cref{sec:tools}, we briefly overview the tools and concepts from information geometry and optimal transport that are directly relevant to the subsequent analysis. In \cref{sec:method}, we summarize the DA-MD approach (with details in \cref{app:MD}) and illustrate how the information-geometric tools and the MD can be naturally combined to reduce predictive uncertainty. \Cref{sec:problem} contains results of our numerical experiments conducted on a Langevin equation with either white or colored noise. Main conclusions drawn from this study are summarized in \cref{sec:concl}.

\section{Preliminaries}
\label{sec:tools}

Let ${}_p \mathcal P (\mathbb R^d)$ denote the probability space of PDFs $f$ on $\mathbb R^d$
%space of probability measures on $\{\mathbb R^d,\mathcal B(\mathbb R^d)\}$ 
with finite $p$th moments, where $p \ge 1$.
%and $\mathcal B(\mathbb R^d)$ is the Borel $\sigma$-algebra on $\mathbb R^d$. 
Our key objective is to minimize loss functions involving PDFs $f$ belonging to ${}_p \mathcal P (\mathbb R^d)$.  In this section, we summarize definitions, tools and theoretical results that will be subsequently used in concert with  DA-MD. 

\paragraph{Measures of discrepancy.}

Alongside classic measures of discrepancy between generic integrable functions $f_1(\mathbf X), f_2(\mathbf X): \mathbb R^d \rightarrow \mathbb R^+$ such as the $L_1$ and the $L_2$ norms, 
\begin{align*}
    d_1(f_1, f_2) := \int_{\mathbb R^d} | f_1(\mathbf X) - f_2(\mathbf X) |\text d \mathbf X  \quad \text{and} \quad  d_2(f_1, f_2) := \left( \int_{\mathbb R^2} | f_1( \mathbf X) - f_2( \mathbf X) |^2 \text d \mathbf X \right)^{1/2}, \label{eq:L2}
\end{align*}
we utilize measures of discrepancy that are tailored to the underlying geometry of the probabilistic space ${}_p \mathcal P (\mathbb R^{d})$.  The KL divergence,
\begin{equation}\label{eq:KL}
    d_{\text{KL}} (f_1, f_2)  := \int_{\mathbb R^d} f_1(\mathbf X) \ln \frac{f_1 (\mathbf X)}{f_2(\mathbf X)} \text d \mathbf X,
\end{equation}
expresses the discrepancy between the PDFs $f_1$ and $f_2$ in terms of relative entropy. Used to quantify how well $f_1 :\mathbb R^d \rightarrow {}_p \mathcal P(\mathbb R^d)$ approximates $f_2 :\mathbb R^d \rightarrow {}_p \mathcal P(\mathbb R^d)$, the KL divergence is not a distance since $d_{\text{KL}} (f_1, f_2) \neq d_{\text{KL}} (f_2, f_1)$.
%confers an invariant geometric structure on the statistical manifold of PDFs $f \in {}_p \mathcal P (\mathbb R^d)$. 
%; specifically, the Hessian operator of the KL divergence turns the probability space ${}_p \mathcal P(\mathbb R)$ into a statistical manifold $\mathcal M$, whose Riemannian metric is the Fisher Information Matrix $\mathbf G_F$ and which enjoys parameterization invariance. 

Another discrepancy measure is the $p$-Wasserstein distance,
\begin{equation} \label{eq:Wp}
    W_p(f_1, f_2) := \left( \underset{\gamma \in \Gamma(f_1, f_2)}{\inf} \int_{\mathbb R^d \times \mathbb R^d} \| \mathbf X - \mathbf Y\|^{p} \gamma(\text d \mathbf X,\text d \mathbf Y) \right)^{1/p}, \quad p \ge 1, 
\end{equation}
where $\Gamma$ is the set of joint probability measures $\gamma$ on $\mathbb R^d \times \mathbb R^d$ whose marginals are univariate probability measures corresponding to $f_1$ and $f_2$. Originating in the field of optimal transport,~\eqref{eq:Wp} 
quantifies the optimal (infimum) cost of shifting the mass distribution of $f_1$ to $f_2$. Such minimum exists and is unique under regularity conditions for the univariate PDFs for $p>1$, i.e., $f$ must be absolutely continuous with respect to the Lebesgue measure~\cite{villani2003topics}.  For $d=1$, \eqref{eq:Wp} reduces to \cite{panaretos2019statistical}
\begin{equation}\label{eq:Wp1d}
    W_p(f_1, f_2) = \| F_1^{-1}(Y) - F_2^{-1}(Y) \|_p = \left( \int_0^1 | F_1^{-1}(Y) - F_2^{-1}(Y)|^p \text d Y \right)^{1/p}, \quad p \ge 1,
\end{equation}
where $F_i(X) = \int_{-\infty}^X f_i(X) \text d X$ with $i=1,2$ is the CDF corresponding to the PDF $f_i(X)$; and $F_i^{-1}(Y)$ is the inverse of $F_i$ defined as $F_i^{-1}(Y) = \inf \{X: F_i(X) \ge Y, Y \in (0,1)\}$. 

Since DA-MD deals with univariate distributions, we are concerned with $d=1$.
 
%
%The  geometric space of PDFs $f \in {}_2 \mathcal P(\mathbb R)$ defined by $W_2$ is characterized by the $L_2$-Wasserstein metric tensor, which endows the manifold of PDFs with an infinite-dimensional Riemannian differential structure. 
%; if finite second moments it works more rigorously. 

\paragraph{Approximation of distributions.}

Various fields of science and engineering---e.g., machine learning~\cite{frogner2015learning,arjovsky2017wasserstein}, estimation theory~\cite{neyman1948consistent}, and optimal transport and control theory~\cite{villani2003topics, esfahani2018data, boso-2020-ambiguity}---deal with a problem of approximating an (empirical) target PDF $\hat f(X)$ with a PDF $f(X;\boldsymbol \varphi) : \mathbb R \rightarrow \mathcal P_{\boldsymbol\varphi} (\mathbb R)$ defined on the parameterized probability space $\mathcal P_{\boldsymbol\varphi}$. The latter consists of PDFs that are uniquely characterized by a set of $N_\text{par}$ parameters $\boldsymbol \varphi \in \Phi \subset \mathbb R^{N_\text{par}}$ with $N_\text{par} \ge 1$. This functional approximation is recast as a problem of finding a parameter set that minimizes a function $\mathcal C(\boldsymbol \varphi)$ depending on a selected measure of discrepancy $D(\boldsymbol \varphi)$ between the target PDF $\hat f(X)$ and its approximation $f(X;\boldsymbol \varphi)$, 
\begin{equation}\label{eq:mindd}
     \underset{\boldsymbol \varphi \in \Phi}{\text{argmin}} \;  \mathcal C (D(\boldsymbol \varphi)), \quad \text{with} \; D(\boldsymbol \varphi) = D(f(X;\boldsymbol \varphi),\hat f(X)), %\quad \text{subject to} \;  f(\boldsymbol \varphi) \in \mathcal P_\varphi.
\end{equation}
with $f(X;\boldsymbol \varphi)$ belonging to $\mathcal P_{\boldsymbol \varphi}$. We assume $\mathcal P_{\boldsymbol \varphi}$ to be a subset of ${}_2{\mathcal P}(\mathbb R)$. The use of the KL and $W_2$ metrics in place of $D$ in \eqref{eq:mindd} introduces known geometries to the statistical manifold of parameterized PDFs, facilitating the deployment of predictable optimization algorithms that exploit this geometric structure. Specifically, one of the geometric properties of the KL divergence is its  parameterization invariance, i.e., the equivalency between computation of the discrepancy $\mathcal C(\boldsymbol \varphi) \equiv D(\boldsymbol \varphi) \equiv d_\text{KL}(f(X;\boldsymbol \varphi), \hat f(X))$ in the PDF space $\mathcal P_\varphi$ and in the parameter space $\Phi$; this property facilitates minimization of the loss function via natural gradient descent~\cite[Sec. 2.1.3]{ollivier2017information}. Moreover, a solution of the minimization problem~\eqref{eq:mindd} with $\mathcal C (\boldsymbol \varphi) \equiv d_{\text{KL}}(f(X;\boldsymbol \varphi),\hat f(X))$ corresponds to the maximum likelihood estimate of the parameters $\boldsymbol \varphi$~\cite{li2017bayesian}. This analogy elucidates the connection between Bayesian inference and information geometry. When $\hat f$ is obtained empirically (e.g., from sampling or repeated experiments), the use of the Wasserstein distance, $\mathcal C (\boldsymbol \varphi) \equiv D^2(\boldsymbol \varphi)/2\equiv   W_2^2(f(X;\boldsymbol \varphi),\hat f(X)) / 2$, is more computationally expedient~\cite{villani2003topics, frogner2015learning, arjovsky2017wasserstein, esfahani2018data}, while possessing geometric properties almost as rigorous as KL \cite{chen2018wasserstein}.   

%\paragraph{Statistical model.} Following \cite{chen2018wasserstein}, we consider a statistical model as defined by the triple $(\Omega, \Phi, f)$: the region where $f$ is non-negative is the support of the distribution $\Omega := \{ X \in \mathbb R | f(X) > 0\}$, where $X$ is the argument of the distribution. We assume compact support $\Omega = \left[ X_{\text{min}}, X_{\text{max}}\right]\subset \mathbb R$, finite-dimensional parameter space $\Phi \subset \mathbb R^n \; \text{with} \; n<+ \infty$, and positive and smooth densities $f$. $f$, and the corresponding CDF $F(X)$, in the following are intended to depend on the parameters $\boldsymbol \varphi \in \Phi$.  

\paragraph{Statistical manifolds.} Let the PDF $f(X; \boldsymbol \varphi)$ be smooth and have a support $\Omega := \{ X \in \mathbb R | f(X) > 0\}$. We assume this support to be compact, $\Omega = \left[ X_{\text{min}}, X_{\text{max}}\right]\subset \mathbb R$, and the dimensionality of the parameter space $\Phi \subset \mathbb R^{N_\text{par}}$ to be finite, $N_\text{par} < + \infty$. An $N_\text{par}$-dimensional manifold is an $N_\text{par}$-dimensional topological space that behaves locally like the Euclidean space $\mathbb R^{N_\text{par}}$. A smooth manifold is equipped with a metric tensor $\mathbf G(\boldsymbol\varphi)$---which facilitates the calculation of distances on the local approximation of the manifold, i.e., the tangent plane---and an affine connection $\nabla_{\boldsymbol\varphi}$---which enables differentiation. The second-order tensor $\mathbf G$ is positive definite and varies smoothly with $\boldsymbol \varphi$. A statistical manifold $\mathcal M$ is a manifold with coordinates $\boldsymbol \varphi = (\varphi^1, \dots, \varphi^{N_\text{par}}) \in \mathbb R^{N_\text{par}}$ where each point represents a PDF with assigned support and defined features. A divergence on the statistical manifold $\mathcal M$ is a non-negative function $D(f(X;\boldsymbol \varphi), f(X;\boldsymbol \varphi')): \mathcal M \times \mathcal M \rightarrow \mathbb R^+$, which is equal to zero if and only if $f(X;\boldsymbol \varphi) \equiv f(X;\boldsymbol \varphi')$ and which can be approximated locally (i.e., when $\boldsymbol \varphi$ and $\boldsymbol \varphi'$ are close) via the components $G_{ij}$ of the second-order tensor $\mathbf G$ as $D(f(\boldsymbol \varphi),f(\boldsymbol \varphi')) = G_{ij}(\boldsymbol \varphi) \Delta \varphi^i \Delta \varphi^j /2 + \mathcal O(| \Delta \boldsymbol \varphi|^3),$
where $\Delta \boldsymbol \varphi = \boldsymbol \varphi - \boldsymbol \varphi'$ and the Einstein summation is implied over the repeated indices $i,j=1,\dots,N_\text{par}$. The tensor $\mathbf G$ defines a Riemannian metric on the statistical manifold $\mathcal M$, and $\mathcal M$ is said to be Riemannian.

% \textbf{The equality between a divergence between distrivutions and a divergence between their parameters, depends on the convexity of the reference function and it is constructed in a specific remapping of the coordinates (couple coord system, for example through the Legendre transformation that has to be one-to-one and diffferentiable).
% The convexity of a function holds upon coordinate transformation if the coordinate transformation is affine.}

%%%%%%%%%%%%%%%%%%%%%%%%%%%%%%%%%%%%%%%%%%%%%%%%%%%%%%%%%%%%%%%%%%%%%%%%%%
\paragraph{Information geometry of statistical manifolds.} %Riemannian structure induced by the KL divergence and natural gradient descent.}

%Minimizing the KL divergence is equivalent to identifying the maximum likelihood estimate of the parameters of the distribution (this is a very old result), with solid convergence proofs. 
If the KL divergence is used to quantify the discrepancy between two PDFs on the manifold $\mathcal M$, then the tensor metric $\mathbf G(\boldsymbol\varphi)$ (a geometric structure) of the space $\mathcal P_{\boldsymbol\varphi}$ of parameterized univariate PDFs $f(X;\boldsymbol \varphi)$ is called Fisher information matrix,
\begin{equation}\label{eq:GF1d}
    \mathbf G_F(\boldsymbol \varphi) = \int_{\Omega} \frac{1}{f( X; \boldsymbol \varphi)} \left( \nabla_{\boldsymbol\varphi} f( X; \boldsymbol \varphi)\right)^\top \nabla_{\boldsymbol\varphi} f( X;\boldsymbol \varphi) \text d  X.
\end{equation}
The resulting statistical manifold $\mathcal M$ is invariant, i.e., for $\boldsymbol\varphi_i \in \Phi$ and $f_i \equiv f(\boldsymbol\varphi_i)$ with $i=1,2$, the divergence $d_{\text{KL}}(f_1, f_2)$ on the manifold $\mathcal M$ equals the distance $|\boldsymbol \varphi_1 - \boldsymbol \varphi_2|$ in the parameter space $\Phi$. This property underpins the Riemannian natural gradient descent (NGD) method (a.k.a. Fisher-Rao gradient descent) for parameter identification \cite[and references therein]{li2018natural}. The method uses
the metric tensor $\mathbf G_F$ as a pre-conditioner for gradient descent algorithms to solve \eqref{eq:mindd} with $\mathcal C \equiv D \equiv d_{\text{KL}}$,
\begin{align}\label{eq:FGD}
    \boldsymbol \varphi_{k+1} = \boldsymbol \varphi_k - \eta \mathbf G_F^{-1}(\boldsymbol \varphi_k) \nabla_{\boldsymbol\varphi} d_{\text{KL}}(f(X;\boldsymbol \varphi),\hat f)|_{\boldsymbol \varphi_k},
\end{align}
where $\eta$ is the descent step and $\mathbf G_F^{-1}$ is the inverse of $\mathbf G_F$. The technique presents strong theoretical analogies with classic filtering techniques (namely Kalman filter and extended Kalman filter)~\cite{ollivier2018online, ollivier2019extended}. In the absence of an analytical expression for $\mathbf G_F$, the matrix can be approximated empirically, although with debatable accuracy \cite{kunstner2019limitations}. 

Geometric structure, including the metric tensor $\mathbf G_W(\boldsymbol \varphi)$, of the finite-dimensional Wasserstein manifolds of Gaussian PDFs was studied in~\cite{takatsu2011wasserstein,malago2018wasserstein}. These results were subsequently generalized to construct $\mathbf G_W(\boldsymbol \varphi)$ for the manifolds $\mathcal M$ of generic discrete~\cite{li2018natural} and continuous~\cite{chen2018wasserstein} distributions. Specifically, when $d=1$, the Wasserstein manifold's metric tensor $\mathbf G_W$ has an explicit form,
\begin{equation}\label{eq:GW1d}
    \mathbf G_W(\boldsymbol \varphi) = \int \frac{1}{f(X;\boldsymbol \varphi)} \left( \nabla_{\boldsymbol\varphi} F(X;\boldsymbol \varphi) \right)^\top \nabla_{\boldsymbol\varphi} F(X;\boldsymbol \varphi)  \text d X.
\end{equation}

Under some mild regularity assumptions, the finite-dimensional Wasserstein manifold $\mathcal M$ in the parameter space $\Phi$ is Riemannian~\cite{chen2018wasserstein}. It introduces an NGD in the space $\Phi$,
\begin{align}\label{eq:WGD}
    \boldsymbol \varphi_{k+1} = \boldsymbol \varphi_k - \eta \mathbf G_W^{-1}(\boldsymbol \varphi_k) \nabla_{\boldsymbol\varphi}  \mathcal C(\boldsymbol \varphi)|_{\boldsymbol \varphi_k}, \quad \text{with} \quad \mathcal C \equiv \frac{1}{2} D^2 \; \text{and} \; D \equiv W_2. %\quad \text{and} \quad W_2(\boldsymbol \varphi) = W_2(f(X;\boldsymbol \varphi), \hat f).
\end{align}
%The numerical cost of~\eqref{eq:WGD} depends on the number of iterations and on the cost of evaluating $\mathbf G_W^{-1}$. 
%\cite{chen2018wasserstein} point out that - unlike the KL divergence - the manifold in the finite-dimensional parameter space is not totally geodesic. % this means that the distance in the parent manifold is not the same as in the child manifold, although they are related. for kl divergence, we have instead invariance. this pull back operation, although not super elegant, is very practical because it allows one to work in a finite dimensional space.

\begin{remark}\label{rem:cost}
Regardless of whether one chooses the KL divergence or the $W_2$ distance, NGD orients the optimization problem~\eqref{eq:WGD} according to the topology of the statistical manifold $\mathcal M$ as expressed by its metric tensor $\mathbf G_i$ ($i = F$ or $W$), thus accelerating the solution. The computational cost of both~\eqref{eq:FGD} and~\eqref{eq:WGD} depends on the overall number of iterations and on the calculation of $\mathbf G_i$ (storage cost $\mathcal O(N_\text{par}^2)$ per iteration) and its inverse $\mathbf G_i^{-1}$ (inversion cost $\mathcal O(N_\text{par}^3)$ per iteration) \cite{ollivier2017information}. Thus, the overall cost of optimization is a trade-off between the number of iterations, arguably reduced on information-geometric grounds, and the cost of inverting the metric tensor $\mathbf G_i$.
\end{remark}

\begin{remark}
%For an infinite-dimensional manifold of smooth PDFs with a compact domain, the $L_2$-Wasserstein distance yields a geometric structure described by a metric tensor $G_W$.  At the same time, 
The finite-dimensional $L_2$-Wasserstein manifold $\mathcal M$ is not exactly geodesic (unless PDFs are Gaussian), and as such the geodesic distance on the manifold is not identical to $W_2$ \cite{chen2018wasserstein}. As demonstrated by \cite[Th. 1 and Prop. 6]{chen2018wasserstein}, the natural gradient trajectory approximates the geodesic distance up to second order information.
\end{remark}

\begin{remark}
A unifying framework connecting the KL and $W_2$ metrics for manifolds of discrete distributions is proposed in \cite{amari2018information, cuturi2013sinkhorn}. 
\end{remark}

%%%%%%%%%%%%%%%%%%%%%%%%%%%%%%%%%%%%%%%%%%%%%%%%%%%%%%%%%
\section{DA-MD with DNN Surrogates}
\label{sec:method}
%%%%%%%%%%%%%%%%%%%%%%%%%%%%%%%%%%%%%%%%%%%%%%%%%%%%%%%%%

Consider a state variable $x(t): \mathbb R^+ \rightarrow \mathbb R$, whose dynamics is governed by a stochastic/random ordinary differential equation
\begin{subequations}\label{eq:ODEphysical}
\begin{align}
    \frac{\text d x(t)}{\text d t} = s(x(t); w(t), \boldsymbol \theta ), \qquad t>0;
\end{align}
subject to a (possibly uncertain, i.e., random) initial condition
\begin{align}
    x(t=0) = x_0, \qquad x_0 \in \mathbb R.
\end{align}
\end{subequations}
The system is driven by the stationary (statistically homogeneous) random process $w(t)$ characterized by a single-point PDF $f_w(W;t)$ and a two-point auto-correlation function $\rho_w(|t_1-t_2|)$; these functions involve meta-parameters ${\boldsymbol\varphi}_w$ such as the mean, variance, and correlation length of $w(t)$. The deterministic function $s(x;\cdot)$, parameterized by a set of $N_\theta$ (possibly uncertain, i.e., random) coefficients ${\boldsymbol\theta}  \in \mathbb R^{N_\theta}$, is such that a solution to~\eqref{eq:ODEphysical} is smooth almost surely in the probability space of both $w(t)$ and, possibly, $\boldsymbol \theta$ and $x_0$. If $\boldsymbol \theta$ and $x_0$ are random, then they are characterized by PDFs $f_{\boldsymbol \theta}(\boldsymbol\Theta)$ and $f_0(X)$, with meta-parameters ${\boldsymbol\varphi}_\theta$ and ${\boldsymbol\varphi}_0$, respectively. In all, the statistics of $x(t)$ depends on the set of $N_\text{par}$ meta-parameters $\boldsymbol \varphi = (\boldsymbol \varphi_w, \boldsymbol\varphi_\theta, \boldsymbol\varphi_0) \in \Phi \subset \mathbb R^{N_\text{par}}$.

In addition to being described by the model~\eqref{eq:ODEphysical}, the system state $x(t)$ is sampled at $N_\text{meas}$ times $t_1,\dots,t_{N_\text{meas}}$. The noisy observations $\hat{\mathbf x} = \{\hat x_1, \dots, \hat x_{N_\text{meas}} \}$ satisfy the data model
\begin{equation}\label{eq:data}
    \hat x_m = x(t_m) + \varepsilon_m, \qquad m=1,\dots,N_\text{meas},
\end{equation}
where the Gaussian measurement errors $\varepsilon_m$ are mutually uncorrelated and have zero mean and variance $\sigma_\varepsilon^2$. 

A goal of data assimilation (DA) is to improve model predictions by augmenting them with observations. Some DA methods yield the ``best'' (i.e., unbiased) prediction and quantify its predictive uncertainty in terms of, respectively, the ensemble mean, $\langle x(t) \rangle$, and the standard deviation, $\sigma_x(t)$, of the state variable $x(t)$. These statistics provide but limited information about $x(t)$, unless its single-point PDF $f(X;t)$ is Gaussian or a known map thereof. Bayesian update and particle filters are examples of DA strategies that overcome this limitation by seeking a solution of~\eqref{eq:ODEphysical} in terms of the PDF $f(X;t)$---or the corresponding CDF $F(X; t) = \mathbb P[x(t) \le X]$---updated with the data $\hat{\mathbf x}$ in~\eqref{eq:data}. Computing such distributions with ensemble methods requires a large number of repeated solves of~\eqref{eq:ODEphysical}, which can be prohibitively expensive.

Data assimilation via DA-MD~\cite{boso-2020-learning} aims to significantly accelerate the computation. Like many other DA strategies, DA-MD comprises two steps: forecast and analysis. The first of these steps relies on the model~\eqref{eq:ODEphysical} and makes a prediction of the system state at time $t$ in terms of $f(X;t)$ or $F(X; t)$. Rather than using, e.g., Monte Carlo simulations, the MD~\cite{tartakovsky2016method} implements this step by deriving a deterministic equation for $f(X;t)$ or $F(X; t)$. Thus, the single-point CDF $F(X; t)$ of the state variable $x(t)$ in~\eqref{eq:ODEphysical} satisfies (sometime approximately) a parabolic PDE (\cref{app:MD})\footnote{For spatially-dependent physical models, space would appear as a coordinate in a CDF or PDF equation \cite{boso-2020-learning}. For systems, the MD would yield a PDF equation for the joint PDF of the interacting system states~\cite{boso-2018-probabilistic,Alawadhi2018Method}.}
\begin{subequations}\label{eq:CDFeq_closed}
\begin{equation}
    \frac{\partial F}{\partial t} + \mathcal U(X,t;\boldsymbol \varphi) \frac{\partial F}{\partial X} = \frac{\partial }{\partial X} \left( \mathcal D(X,t;\boldsymbol \varphi) \frac{\partial F}{\partial X} \right), \qquad t>0, \quad X \in \Omega = [X_{\text{min}}, X_{\text{max}}],
\end{equation}
subject to initial and boundary conditions
\begin{align}\label{eq:bc}
F(X;0) = F_0(X), \qquad F(X_{\text{min}},t) = 0, \qquad F(X_{\text{max}},t) = 1.
\end{align}
\end{subequations}
The drift velocity, $\mathcal U(X,t;\boldsymbol \varphi) : \Omega \times \mathbb R^+ \rightarrow \mathbb R$, and the diffusion coefficient, $\mathcal D(X;t,\boldsymbol \varphi) : \Omega \times \mathbb R^+ \rightarrow \mathbb R^+$, are smooth functions of their arguments, which involve a set of the meta-parameters $\boldsymbol\varphi$. The functional forms  $\mathcal U$ and $\mathcal D$ depend on that of $s(x;\cdot)$, on the statistical characterization of the random parameters epitomized by the statistical parameters ${\boldsymbol\varphi}$ of their distributions, and on the degree of approximation introduced by the closure strategy. If the initial state of the system, $x_0$, is known with certainty, then its CDF $F_0(X)$ is the Heaviside step function, $F_0(X) = \mathcal H(X-x_0)$. 

\begin{remark}
The CDF equation~\eqref{eq:CDFeq_closed} maps the meta-parameters ${\boldsymbol \varphi}$ onto $F(X; t, \boldsymbol \varphi)$, the CDF of the system state $x(t)$. In other words, a point $\boldsymbol \varphi \in \Phi \subset \mathbb R^{N_\text{par}}$ can be thought of as a coordinate on the statistical manifold $\mathcal M$ of the CDF $F(X; t, \boldsymbol \varphi)$ at time $t$. At any time $t'$, a solution to \eqref{eq:CDFeq_closed} provides an estimate of the CDF $F(X; t',\boldsymbol \varphi)$ dependent on the current characterization of the random inputs expressed by $\boldsymbol \varphi$. Equivalently, points $\widetilde{\boldsymbol \varphi} = \{t,{\boldsymbol \varphi} \}$ define a dynamic statistical manifold $\mathcal M_t$ of the CDF $F(X; \widetilde{\boldsymbol \varphi})$.
\end{remark}

The second step of DA-MD, analysis via Bayesian update, is performed sequentially for each of the $N_\text{meas}$ measurements $\hat x_m$ in~\eqref{eq:data}.  At $m$th assimilation step, the updated meta-parameters $\boldsymbol \varphi^{(m)}$ are computed by solving the minimization problem~\eqref{eq:mindd} for the discrepancy $D$ between the CDF $F(X;t_m,\boldsymbol \varphi)$ predicted by the model~\eqref{eq:CDFeq_closed} and the observational CDF obtained with Bayes' rule,
\begin{align}\label{eq:observational}
    \hat F(X;t_m) = \int_{X_{\text{min}}}^X \hat f(X;t_m)\text d X
    \quad \text{with} \quad
    \hat f(X;t_m) = \frac{f_L(\hat x_m | x(t_m) = X) f(X;t_m,\boldsymbol \varphi^{(m-1)})}{\int_\Omega f_L(\hat x_m | x(t_m) = X) f(X;t_m,\boldsymbol \varphi^{(m-1)}) \text d X}.
\end{align}
Here the likelihood function $f_L(\hat x_m | x(t_m) = X)$ specifies a data model; and the PDF $f(X;t_m,\boldsymbol \varphi^{(m-1)})$, computed by solving the CDF equation~\eqref{eq:CDFeq_closed} with the parameter set $\boldsymbol \varphi^{(m-1)}$ from the previous assimilation step, serves as a prior. In~\cite{boso-2020-learning}, the discrepancy $D$ was expressed in terms of the $L_2$ norm; a consequence of this choice was significant computational cost of solving the minimization problem~\eqref{eq:mindd}. A main innovation of this study is to exploit the geometric structure of the statistical manifolds in the parameter space $\Phi$ by using either the KL divergence~\eqref{eq:KL} or the Wasserstein distance~\eqref{eq:Wp1d} at each assimilation time. This enables us to solve~\eqref{eq:mindd} via NGD, which we henceforth refer to as NGD-KL and NGD-W$_2$ depending on which metric is used. The update of the meta-parameters $\boldsymbol \varphi$ is done using NGD-KL~\eqref{eq:FGD} or NGD-W$_2$~\eqref{eq:WGD}, taking advantage of the explicit formulations for the manifold's metric tensors $\mathbf G_F$ in~\eqref{eq:GF1d} and $\mathbf G_W$ in~\eqref{eq:GW1d}. 
%
%Minimization 
%\begin{equation}\label{eq:minddd}
%    \boldsymbol \varphi^{(m)} = \underset{\boldsymbol \varphi}{\text{argmin}}  \quad D(\boldsymbol \varphi), \quad \boldsymbol \varphi \in \Phi
%\end{equation}
%is thus performed by choosing either $D(\boldsymbol \varphi ) =  d_{\text{KL}}(f(X;t_m,\boldsymbol \varphi),\hat f(X;t_m)) =  d_{\text{KL}} (\boldsymbol \varphi )$ or $D(\boldsymbol \varphi) = \frac{1}{2} W_2(f(X;t_m,\boldsymbol \varphi),\hat f(X;t_m))^2 = \frac{1}{2} W_2(\boldsymbol \varphi)^2 $ (referred to as KL/W$_2$ loss function subject to KL/W$_2$ optimization, respectively, in the following), and the corresponding Natural Gradient Descent (NGD) algorithm (termed NGD-KL and NGD-W$_2$, respectively). 

\begin{remark}
The analysis step of DA-MD is performed on univariate (one-point) distributions ($d=1$) regardless of the size of the physical parameter and meta-parameter sets, $N_\theta$ and $N_\text{par}$. That drastically reduces (to one) the dimensionality of the update effort in classical Bayesian DA.  Moreover, availability of a CDF/PDF equation removes the need for Gaussianity and linearity assumptions on the physical model and its random parameters. The CDF/PDF equation is assumed to be valid throughout the assimilation process. 
\end{remark}

\begin{remark}
Parameter update via discrepancy minimization places DA-MD in the company of many machine-learning and optimal-transport techniques (see the references above). Unlike these methods, DA-MD uses CDF or PDF equations and their parameters to define the parameter space for a statistical manifold $\boldsymbol \varphi$, such that the discrepancy minimization is constrained by these PDEs. Learning occurs on the statistical manifold defined by $\boldsymbol \varphi$ and proceeds by sequential updates of these meta-parameters.
\end{remark}

%The reliance on a measure of discrepancy associates DA-MD with machine learning and optimal transport theory, with the distinction that the parameter space for the statistical manifold is defined via the CDF equation and its parameters $\boldsymbol \varphi$, such that minimization is constrained by the CDF equation. Learning occurs on the statistical manifold defined by $\boldsymbol \varphi$ and proceeds by sequential updates of the parameters $\boldsymbol \varphi^{(m)}$. The dimensionality of $\boldsymbol \varphi$, $n$, is finite and typically relatively small.
%The dynamic dimension (subsequent measurements in time) is handled sequentially, updating each measurement separately. 

%Since we are always dealing with uni-variate distributions, it is convenient to use either the KL divergence or the L2 Wasserstein metric. Both are associated with geometric shapes of the manifold, which can be used for preconditioning of the gradient descent, transforming into a natural gradient descent.

 \begin{figure}[htbp]
    \centering
    \includegraphics[width=.5\textwidth]{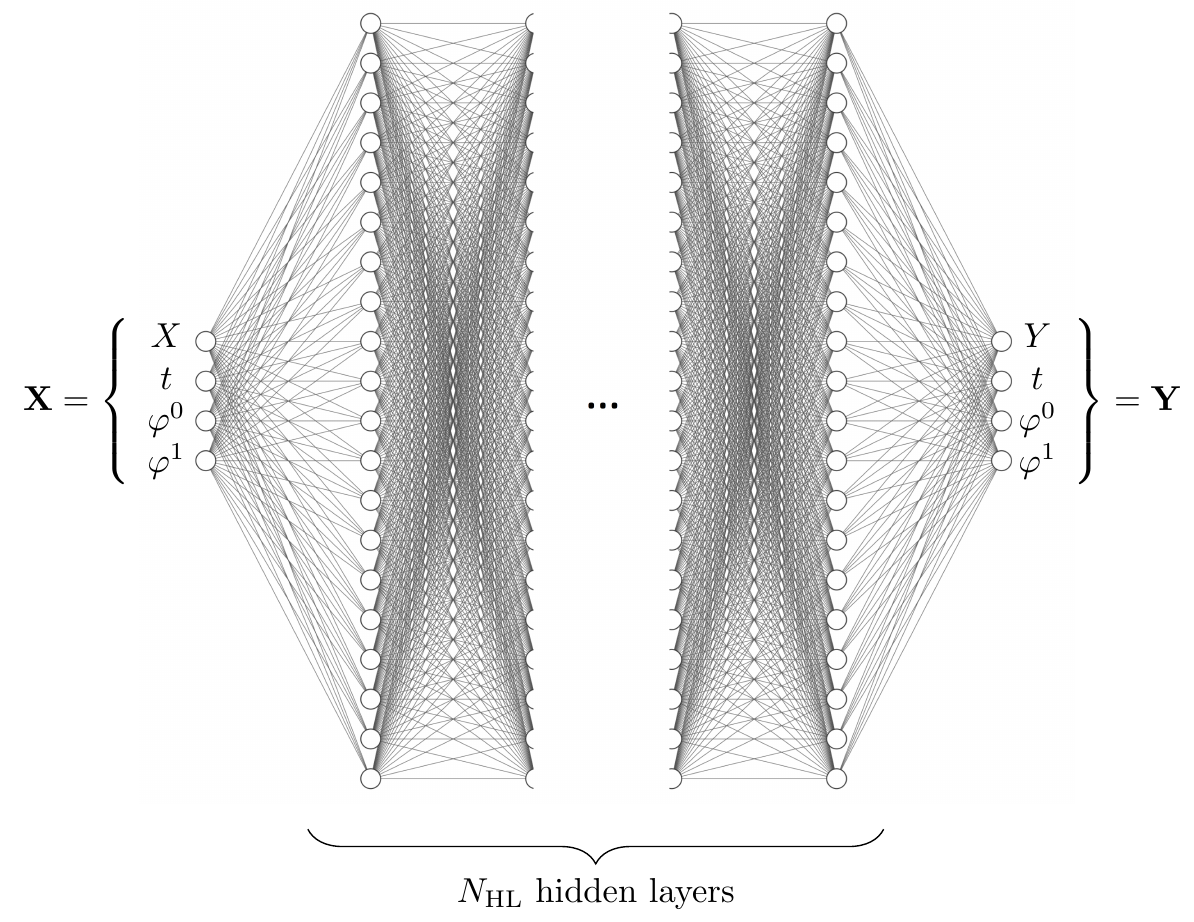}
    \caption{Fully-connected DNN used to approximate a solution of the CDF equation~\eqref{eq:CDFeq_closed}, with $N_{\text{HL}}$ hidden layers and $N_{\text{nphl}}$ nodes per hidden layer.  Inputs are $\mathbf X \{ X = F^{-1}(Y; t, \boldsymbol \varphi), t, \boldsymbol \varphi, \}$, and outputs are $\mathbf Y = \{ Y = F(X; t, \boldsymbol \varphi), t, \boldsymbol \varphi \}$. This illustration has $N_\text{par} = 2$, with $\boldsymbol \varphi = \{ \varphi^0, \varphi^1 \}$. The DNN parameters, in both the illustration and our numerical examples, are set to $N_{\text{HL} } = 7$, $N_{\text{nphl}} = 20$.}
    \label{fig:fcnn}
\end{figure}

\paragraph{Loss function minimization.} 

We use a surrogate model to accelerate the calculation of the discrepancies  $d_\text{KL}(f(X;\boldsymbol \varphi),\hat f)$ or $W_2(f(X;\boldsymbol \varphi),\hat f)$, their respective  gradients $\nabla_{\boldsymbol\varphi} d_\text{KL}$ or $\nabla_{\boldsymbol\varphi} W_2$, and the preconditioning tensor metrics $\mathbf G_F$ or $\mathbf G_W$. Specifically, a fully-connected deep neural network (DNN), whose architecture is  %\cite{raissi2019physics} with a specific architecture (
illustrated in \cref{fig:fcnn}, is used to approximate the solution of the CDF equation given the set of inputs $\mathbf X = \{X,t,\boldsymbol \varphi\}$. The number of outputs in this DNN equals the number of inputs, $\mathbf Y = \{Y_j: j = 1, \dots, N_\text{par} + 2 \} = \{ Y = F(X;t,\boldsymbol \varphi),t,\boldsymbol \varphi\}$, such that $\text{dim}({\mathbf X}) = \text{dim}( \mathbf Y) = N_\text{par} + 2$. We require the resulting vector function $\mathbf Y = \mathbf F(\mathbf X)$ to be one-to-one except at singularity points, and its derivative to be an invertible linear map in local, identifiable regions. These requirements fulfill the hypotheses of the Inverse Function Theorems \cite[Th. 1-2 in sec. 3.2]{guzman2012derivatives} for vector functions.  
Under these conditions, the vector function $\mathbf Y = \mathbf F(\mathbf X)$  is invertible, its inverse is differentiable, and the derivative of the inverse is equal to the inverse of the derivative \cite[Th. 3 in sec. 3.2]{guzman2012derivatives}. %\marginFB{Add more details about these theorems and their hypotheses.}
Automatic differentiation is employed both to verify the inversion theorem hypotheses and to calculate the terms appearing in the minimization algorithms. This is especially useful, since NGD-KL utilizes the derivatives of the forward pass, whereas NGD-W$_2$ requires the derivatives of the inverse function.  A differentiable DNN allows accurate calculation of the metric tensors for both geometries, eliminating potential problems related to their empirical approximation.

The DNN is trained on a data set consisting of $N_{\text{ts}}$ pairs $(\mathbf X_{\text{ts}}^i, \mathbf Y_{\text{ts}}^i)$, for $i=1,\dots,N_{\text{ts}}$. This training set is generated by solving the CDF equation~\eqref{eq:CDFeq_closed} for $N_{\text{ts}}$ combinations of meta-parameters $\boldsymbol \varphi$, i.e., at points $\boldsymbol\varphi_i \in \Phi$ with $i=1,\dots,N_{\text{ts}}$.\footnote{For each $i$, the data pairs $(\mathbf X_{\text{ts}}^i, \mathbf Y_{\text{ts}}^i)$ are extracted from these solutions at  regularly-spaced time intervals and at spatial locations (in the $X$ direction) refined with a cosine mapping around a solution of \eqref{eq:ODEphysical} with mean parameters.} 
The DNN training is accomplished by solving an optimization problem~\cite{raissi2019physics}, 
\begin{subequations}\label{eq:losssurrogatemodel}
\begin{align}
    \underset{\mathbf w,\mathbf b}{\text{argmin}}  ( \text{MSE}_{\text{ts}} + \text{MSE}_{\text R} + \text{MSE}_{\text{aux}} + \text{SMR}),
\end{align}
with respect to the weights and biases of the DNN, $\mathbf w$ and $\mathbf b$, respectively. Here,
\begin{align}
& \text{MSE}_{\text{ts}} = \sum_{j=1}^{n + 2}  \lambda_j \frac{1}{N_{\text{ts}}} \sum_{i=1}^{N_{\text{ts}}}  |  Y_j(\mathbf X_{\text{ts}}^i) -  Y_{j,\text{ts}}^i |^2, \quad \lambda_j = \left( \max  Y_{j,\text{ts}} \right)^{-1}  \\
    & \text{MSE}_{\text{R}} = \frac{1}{N_{\text{R}}} \sum_{i=1}^{N_{\text{R}}} | R(\mathbf X_{\text{R}}^i) |^2  \label{eq:resR}\\
    & \text{MSE}_{\text{aux}} = \frac{1}{N_{\text{aux}}} \sum_{i=1}^{N_{\text{aux}}} |  Y(\mathbf X_{\text{aux}}^i) -  Y_{\text{aux}}^i |^2  \\
    & \text{SMR} = \left( \max \left | \frac{\partial Y}{\partial X} (\mathbf X_{\text{SMR}}) \right | \right)^{-1} \sum_{i=1}^{N_{\text{SMR}}} \max \left( 0, - \frac{\partial Y}{\partial X}(\mathbf X_{\text{SMR}}^i) \right),  \label{eq:SMR}
\end{align}
\end{subequations}
and $\mathbf Y(\mathbf X^i)$ represents the $N_\text{par}+2$ outputs of the DNN with inputs $\mathbf X^i$. The mean square errors $\text{MSE}_{\text{R}}$ and $\text{MSE}_{\text{aux}}$ enforce the fulfillment of the CDF equation and its initial/boundary conditions at collocation points $\{ \mathbf X_{\text R}^i \}_{i=1}^{N_{\text R}}$ and $\{ \mathbf X_{\text{aux}}^i \}_{i=1}^{N_{\text {aux}}}$, respectively.\footnote{We select a regularly spaced set of points for the enforcement of \eqref{eq:resR} in all but the $X$ direction, wherein points are refined around the solution of \eqref{eq:ODEphysical} with mean parameters; $N_{\text{aux}}$ points are regularly spaced in all directions.} The residual is defined as
\begin{equation}\label{eq:res}
    R(\mathbf X_{\text R}^i) = \frac{\partial Y( \mathbf X_{\text R}^i ) }{\partial t} + 
    \left( \mathcal U(\mathbf X_{\text R}^i ) - \frac{\partial \mathcal D } {\partial X} (\mathbf X_{\text R}^i ) \right) 
    \frac{\partial Y }{\partial X} (\mathbf X_{\text R}^i ) - \mathcal D ( \mathbf X_{\text R}^i ) \frac{\partial^2 Y}{\partial X^2} (
    \mathbf X_{\text R}^i ),
\end{equation}
and $Y_{\text{aux}}^i$ represent the auxiliary conditions for the CDF equation at  points $\mathbf X_{\text{aux}}^i$, which represents initial or boundary conditions~\eqref{eq:bc}. The term SMR is a soft constraint \cite{gupta2019incorporate} that regularizes the DNN by enforcing monotonicity of the output $Y = F(X;t,\boldsymbol \varphi)$ along the $X$ direction at points $\{ 
\mathbf X_{\text{SMR}}^i \}_{i=1}^{N_{\text{SMR}}} = \{ \{ \mathbf X_{\text {ts}}^i \}_{i=1}^{N_{\text R}},
\{ \mathbf X_{\text R}^i \}_{i=1}^{N_{\text R}}, \{ \mathbf X_{\text{aux}}^i \}_{i=1}^{N_{\text {aux}}} \}
$. 
The physics-aware component of \eqref{eq:losssurrogatemodel}, $\text{MSE}_{\text R} + \text{MSE}_{\text{aux}}$, makes training less data-intensive and increases confidence in the predictions of the DNN outside the training range (but within the residual points range). %It is good practice to maintain a balance between $N_{\text{TS}}$ and $N_{\text R}$, the size of the training set and residual points ensemble, respectively. %; increasing both does not affect unreasonably the training time. 

\section{Numerical Experiments}
\label{sec:problem}

In this section, we apply the information-theoretic DA strategy introduced above to three problems described by~\eqref{eq:ODEphysical}. Section~\ref{sec:problem1} contains an example of deterministic nonlinear dynamics starting from a random initial condition; this setting provides an ideal testbed for the information-geometric analysis by virtue of lending itself to analytical treatment. \Cref{sec:problem2} deals with a Langevin equation with white noise $w(t)$, a problem for which the CDF equation~\eqref{eq:CDFeq_closed} is exact. In other words, the forecast component of DA-MD is exact, whereas the analysis step introduces an approximation. In \cref{sec:problem3}, we consider a Langevin equation with colored noise $w(t)$ that is modeled as an Ornstein-Uhlenbeck process; the derivation of the CDF equation~\eqref{eq:CDFeq_closed} requires a closure approximation. In this case, the performance of DA-MD depends also on the accuracy and robustness of the CDF equation as forecasting tool.

%We consider \eqref{eq:ODEphysical} with deterministic physical parameters $\{ \boldsymbol \theta_{s_d}, \boldsymbol \theta_{s_w} \}$ and additive random forcing
%\begin{align}\label{eq:ODE}
%    \frac{\partial x}{\partial t} = s_d(x,t;\boldsymbol \theta_{s_d}) + s_w(x,t;\boldsymbol \theta_{s_w}) w(t) , \quad \text{subject to} \quad  x(t=0) = x_0, % v(x(t),t; \boldsymbol \theta_v)
%\end{align}
%and three scenarios: i) deterministic dynamics ($s_w = 0$) and random initial condition; ii) Langevin equation with deterministic initial condition, where $w$ is standard Gaussian white noise (delta-correlated); iii) Langevin equation with colored noise, where $w(t)$ is modeled as a Ornstein Uhlenbeck process and $x_0$ is deterministic.

In all cases, one realization ($x_0^\star$ or $\boldsymbol\theta^\star$) of the relevant random parameters, $x_0$ or $\boldsymbol\theta$, represents ground truth. Statistical models for these parameters are chosen such that the state variable $x(t)$ has a compact support $\Omega \subset \mathbb R^+$. This ensures that the information geometry induced by the $W_2$ divergence is rigorously defined. %We seek to utilize $N_m$ observations of $x(t)$, $\hat {\mathbf x} = \{ x_1, \dots, x_{N_m}\}$, collected at times $t_1, \dots, t_{N_m}$, to refine the knowledge of the meta-parameters $\boldsymbol \varphi$ in both scenarios, thus reducing predictive uncertainty. 
The $N_\text{meas}$ observations $\hat{\mathbf x}$ are taken at regular time intervals, with the time step $\Delta t = t_{N_\text{meas}}/(N_\text{meas}
+1)$. They are generated by adding zero-mean Gaussian noise with standard deviation $\sigma_\varepsilon$ to the solution of~\eqref{eq:ODEphysical} with $x_0^\star$ or $\boldsymbol\theta^\star$ (i.e., the synthetic truth). This procedure results in the Gaussian likelihood function $f_L$, although other choices are possible. While not investigated here, data models constructed on repeated observations of the same phenomenon might be more suitable for processes that are inherently random like those described by Langevin equations.

For the Langevin scenarios in sections~\ref{sec:problem2} and~\ref{sec:problem3}, we employ the \texttt{JITCSDE} Python module \cite{jitcxde} to solve the stochastic ordinary differential equation~\eqref{eq:ODEphysical}. The corresponding CDF equations~\eqref{eq:CDFeq_closed} are solved with a finite volumes (FV) scheme, implemented using the \texttt{Fipy} library \cite{FiPy:2009}, to provide a training set for the surrogate model. DNN is trained by employing \texttt{Tensorflow}; optimization in~\eqref{eq:losssurrogatemodel} is performed using L-BFGS-B method \cite{liu1989limited}, with a random initialization of $\mathbf w$ and $\mathbf b$;  and the network topology is shown in \cref{fig:fcnn}. Automatic differentiation is used to compute both the derivatives in the residual $R$ in~\eqref{eq:res} and the PDF from CDF.
Minimization of the KL and $W_2$ discrepancies is performed using both standard gradient descent (GD) and NGD. In the case of NGD, convergence is accelerated by the use of the pre-conditioners $\mathbf G_F$ and $\mathbf G_W$ in~\eqref{eq:FGD} and~\eqref{eq:WGD}. For each direction established by the gradient of the loss function (adjusted by the pre-conditioners when NGD is used) we employ the \texttt{Scipy} library's implementation of step calculation \cite[Sec. 5.2]{wright1999numerical}. A convergence criterion for NGD in~\eqref{eq:FGD} and~\eqref{eq:WGD} is defined by 
%\begin{equation}
   $| \nabla_\varphi D | \le \epsilon$.
%\end{equation}
Because of the different order of magnitude of the KL and $W_2$ discrepancies $D$, the convergence threshold $\epsilon$ is discrepancy-specific; we select a KL-based minimization threshold, $\epsilon_{\text{KL}}$, and assign the threshold for W$_2$, $\epsilon_{\text W_2}$, such that $\epsilon_{\text W_2} / \mathcal C(\text W_2(f(X;t_1,\boldsymbol \varphi^{(0)}); \hat f(X;t_1))) = \epsilon_{\text{KL}} / \mathcal C(d_{\text{KL}}(f(X;t_1,\boldsymbol \varphi^{(0)}); \hat f(X;t_1)))$.
%, and assigned so that the NGD computational time over the assimilation window is the same for both divergences.
% \textcolor{green}{
% Like before, given the different order of magnitude of the values of both the loss function and its gradients for the two loss functions, we compare convergence results by imposing a convergence threshold $\epsilon$ (in terms of the gradient of the loss function 
% $| \nabla_\varphi D(\boldsymbol \varphi) |$) such that the DG computational time is the same (within $5 \%$) 
% for both divergences. For the specific settings, this yielded 
% $\epsilon_{\text{KL}} = 10^{-4}$ and $\epsilon_{{W_2}} = 10^{-3}$ for the KL and the W$_2$ loss functions, respectively.}
%

\subsection{Deterministic dynamics with random initial state}
\label{sec:problem1}

The dynamics of state variable $x(t)$ is described by
\begin{align}\label{eq:ODE}
    \frac{\text d x}{\text d t} = - 2 x^2, \qquad x(0) = x_0, 
\end{align}
%The first case lends itself to analytical treatment, and it serves as a testbed for the proposed procedure; we employ $s_d = - a \alpha x^\alpha$, with chosen $a \in \mathbb R^+$ and 
%where $\alpha \in \mathbb N^+$. %regulating the linearity (with $\alpha=1$) or non-linearity (with $\alpha>1$) of the right-hand-side.
The random initial state $x_0$ has compact support $\Omega_0 \subset \mathbb R^+$, which ensures that $x(t)$ has a compact support $\Omega \subset \mathbb R^+$. % and, hence, that the information geometry induced by the $W_2$ divergence is rigorously defined. \textbf{[Supp?]} 
To be specific, and without loss of generality, we take the CDF of $x_0$, $F_0(X;\boldsymbol\varphi_0)$, to be Gaussian, with assigned prior mean ($\mu_0^{(0)}$) and standard deviation ($\sigma_0^{(0)}$) acting as the sole meta-parameters for the model, i.e., $\boldsymbol \varphi_0^{(0)} = \{ \mu_0^{(0)}, \sigma_0^{(0)}\} = \boldsymbol \varphi^{(0)}$. %(Any other choice of distribution could be made.)

For this problem, the general CDF equation~\eqref{eq:CDFeq_closed} is exact, reduces to (\cref{sec:case1})
\begin{align} \label{eq:CDF1}
     \frac{\partial F}{\partial t} - 2 X^2 \frac{\partial F}{\partial X} = 0, \qquad
     F(X;t=0) = F_0(X);
\end{align}
%The exact CDF equation - as detailed in \cref{app:MD} - corresponds to \eqref{eq:CDFeq_closed} with $\mathcal U = s(X,t; a) $ and $\mathcal D = 0$. 
and has an analytical solution $F(X;t, \boldsymbol\varphi)$ and the corresponding analytical expression for the PDF $f(X;t, \boldsymbol\varphi) = \text d F / \text dX$. As a consequence, there is no need for a surrogate model of the solution to this CDF equation. 
The data assimilation problem has a computable Bayesian solution
\begin{equation}\label{eq:Bayesparameter}
    f_0 (X_0 | \hat {\mathbf x}) = \frac{f_L(\hat{\mathbf x} | x(t_{1:N_\text{meas}},X_0))  f_0 (X_0)}{\int_{\Omega} f_L(\hat {\mathbf x} | x(t_{1:N_\text{meas}},X_0)) \hat f_0 (X_0) \text d X_0} = \frac{ \prod_{m=1}^{N_\text{meas}} f_L(\hat {x}_m | x(t_m,X_0) ) f_0 (X_0)} {\int \prod_{m=1}^{N_\text{meas}} f_L(\hat {x}_m | x(t_m,X_0) ) f_0 (X_0) \text d X_0},
\end{equation}
where the prior PDF $f_0(X)$ is computed as $f_0 = \text d F_0 / \text d X$, and the i.i.d. measurements $\hat x_m$ are assigned the likelihood function 
\[f_L(\hat{\mathbf x}| x(t_{1:N_\text{meas}}, X_0)) = \prod_{m=1}^{N_\text{meas}} f_L(x_m | x(t_m,X_0)).
\]

We use the exact Bayesian posterior~\eqref{eq:Bayesparameter} to gauge the accuracy of the sequential Bayesian update of the meta-parameters $\boldsymbol \varphi = \boldsymbol \varphi_0$ via GD for~\eqref{eq:mindd} with $\mathcal C \equiv d_\text{KL}$ or  $W_2^2/2$,  NGD-KL~\eqref{eq:FGD}, and NGD-W$_2$~\eqref{eq:WGD}. %Per the DA-MD procedure, each observation $\hat x_m$ is assimilated sequentially to update the meta-parameters for the CDF equation, $\boldsymbol \varphi = \boldsymbol \varphi_0$. 
Assigned meta-parameters $\boldsymbol \varphi_0$ uniquely identify a distribution for the state $x(t)$ through the (analytical) solution to the CDF equation \eqref{eq:CDF1}. A forecast PDF at the measurement time $t_m$, and the corresponding observational PDF $f(X;t_m)$ is obtained via Bayes'rule 
\begin{equation}\label{eq:Bayesstate}
    \hat f(X;t_m | x_m) = \frac{f_L(\hat{x}_m | x(t_m) = X ) f (X;t_m, \boldsymbol \varphi^{(m-1)})}{\int_{\Omega} f_L(\hat{x}_m | x(t_m) = X ) f (X;t_m, \boldsymbol \varphi^{(m-1)}) \text d X},
\end{equation}
in which the priors $f (X;t_m, \boldsymbol \varphi^{(m-1)}$ are computed analytically. %and predictions $F(X;t, \boldsymbol \varphi)$ to be used in \eqref{eq:mindd} are obtained as analytical solutions of \eqref{eq:CDF1} with Gaussian $F_0(X;\boldsymbol \varphi_0)$. 
The availability of analytical expressions for $F(X;t)$ and $f(X;t)$ facilitates the (semi-)analytical computation of both the metric tensors $\mathbf G_F$ and $\mathbf G_W$ in~\eqref{eq:GF1d} and~\eqref{eq:GW1d}, and the the gradient of the discrepancy, $\nabla_\varphi \mathcal D$, for the KL an $W_2$ measures.
%
%the minimization of the loss function at each step $m$, $D^{(m)} = D(f^{(m)},\hat f^{(m)})$ in \eqref{eq:mindd}, for both KL and $\text W_2$ divergence, since both the gradient of the loss function, $\nabla_\varphi \mathcal D$, and the metric tensors ($\mathbf G_F$ and $\mathbf G_W$ for KL and W$_2$, respectively) can be computed (semi-)analytically. 
The integrals in the metric tensors, the discrepancy gradient, and the normalization constant in \eqref{eq:Bayesstate}, are computed via numerical quadrature from the Fortran library \texttt{QUADPACK}.

\begin{figure}[htbp]
\begin{center}
\begin{subfigure}{.49\textwidth}
    \includegraphics[width=\textwidth]{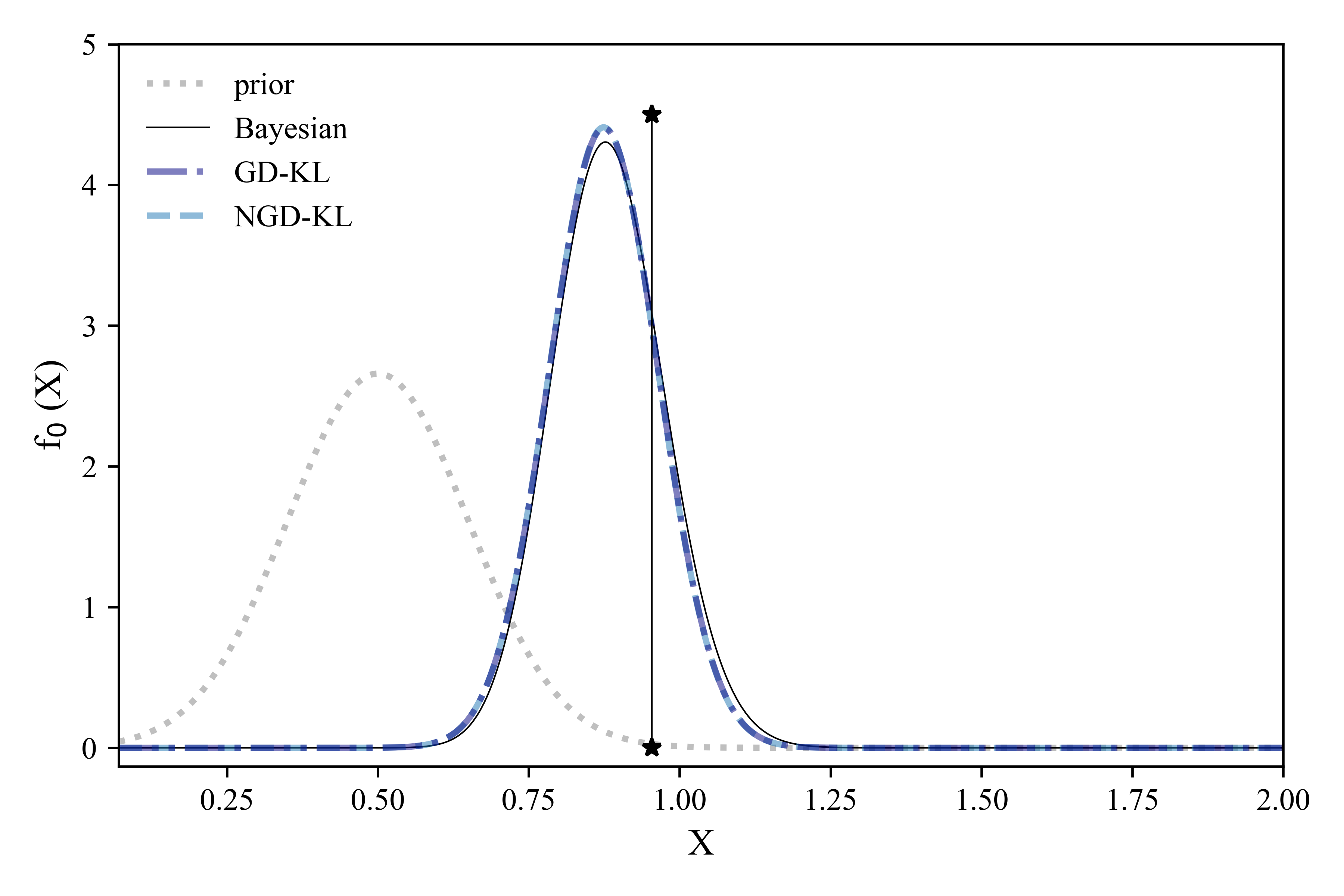}
\end{subfigure}
\begin{subfigure}{.49\textwidth}
    \includegraphics[width=\textwidth]{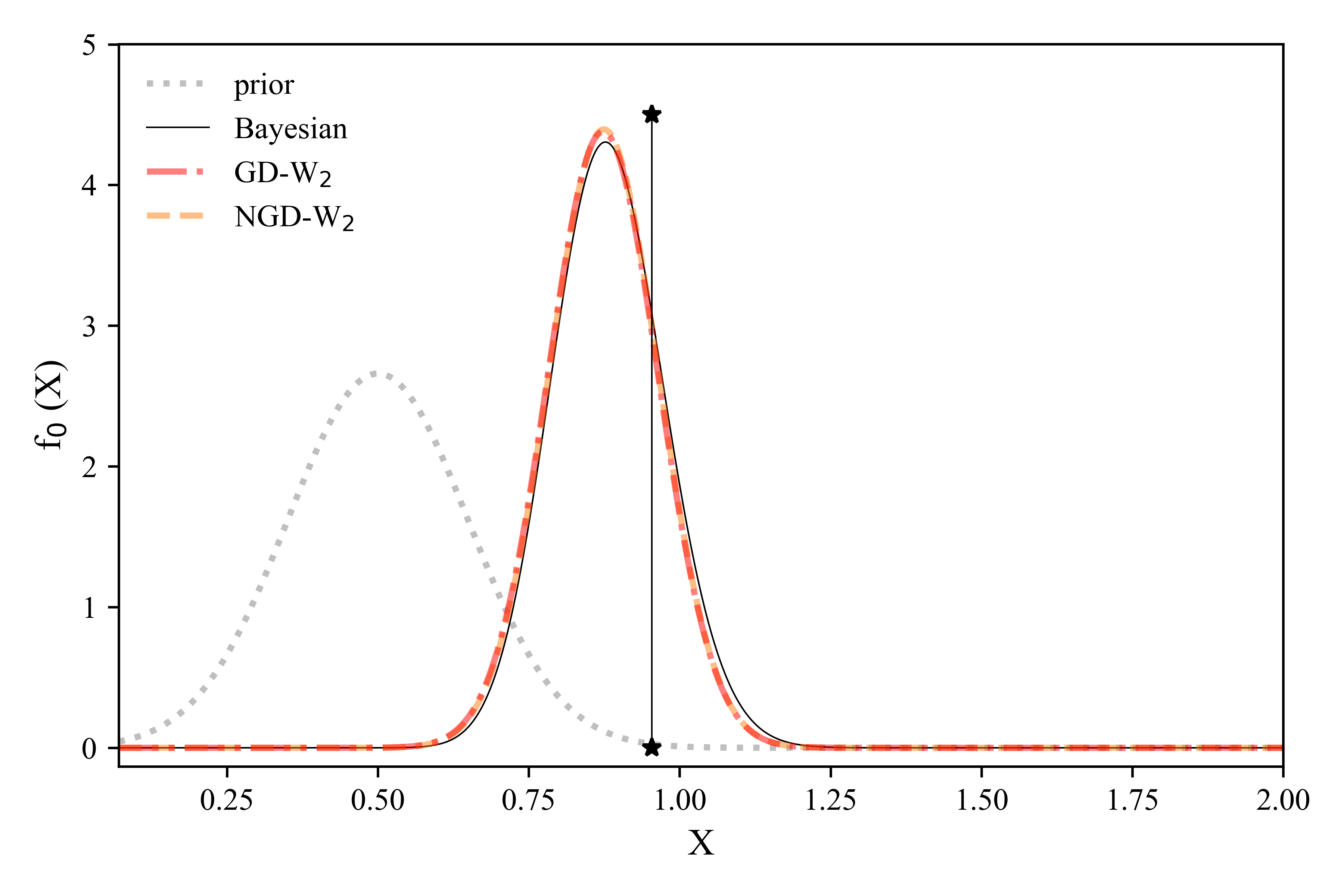}
\end{subfigure}
\end{center}
\caption{Prior and posterior PDFs of the initial state $x_0$ after assimilation of $N_\text{meas} = 10$ observations, computed with the information-geometric optimization strategies NGD-KL (left) and NGD-W$_2$ (right). Also plotted are the analytically derived Bayesian posterior (solid line) and the posteriors obtained with GD-KL (dash-dotted line on the left) and GD-W$_2$ (dash-dotted line on the left). The ground truth, $x_0^\star =0.954$, is indicated by the starred vertical line. Simulation parameters are set to $\sigma_\varepsilon=0.1$,  $t_{N_\text{meas}} = 2$, $\boldsymbol \varphi^{(0)} = (0.5,0.15)$, $\epsilon_{\text{KL}} = 10^{-3}$.}
\label{fig:pinitKL}
\end{figure}

Figure~\ref{fig:pinitKL} exhibits prior and posterior PDFs of the random initial state $x_0$, obtained alternatively with the four DA-MD implementations---GD for~\eqref{eq:mindd} with $D \equiv d_\text{KL}$ or $W_2$, NGD-KL~\eqref{eq:FGD}, and NGD-W$_2$~\eqref{eq:WGD}---and with the analytical Bayesian update~\eqref{eq:Bayesparameter}. %Results for both KL and W$_2$ loss functions (left and right panels, respectively) are shown. %for the non-linear case ($\alpha=2$), although similar plots can be obtained for the linear case. 
The information-geometric optimization strategies NGD-KL and NGD-W$_2$ have comparable performance, both reproducing accurately the exact Bayesian posterior and having negligible difference in the identified meta-parameters $\boldsymbol \varphi$ at the end of the assimilation window. After assimilation of $N_\text{meas} = 10$ measurements, the unknown ground truth $x_0^\star = 0.954$ is approximated by the mean of the posterior PDF, $\mu_0^{(N_\text{meas})} \equiv \varphi_1^{(N_\text{meas})}$; the standard deviation of this PDF, $\sigma_0^{(N_\text{meas})} \equiv \varphi_2^{(N_\text{meas})}$, provides a measure of predictive uncertainty. These statistics are $\boldsymbol \varphi^{(N_\text{meas})} = \{ 0.88,0.09\}$ for all optimization algorithms.

\begin{figure}[htbp]
\begin{center}
\begin{subfigure}{.49\textwidth}
    \includegraphics[width=\textwidth]{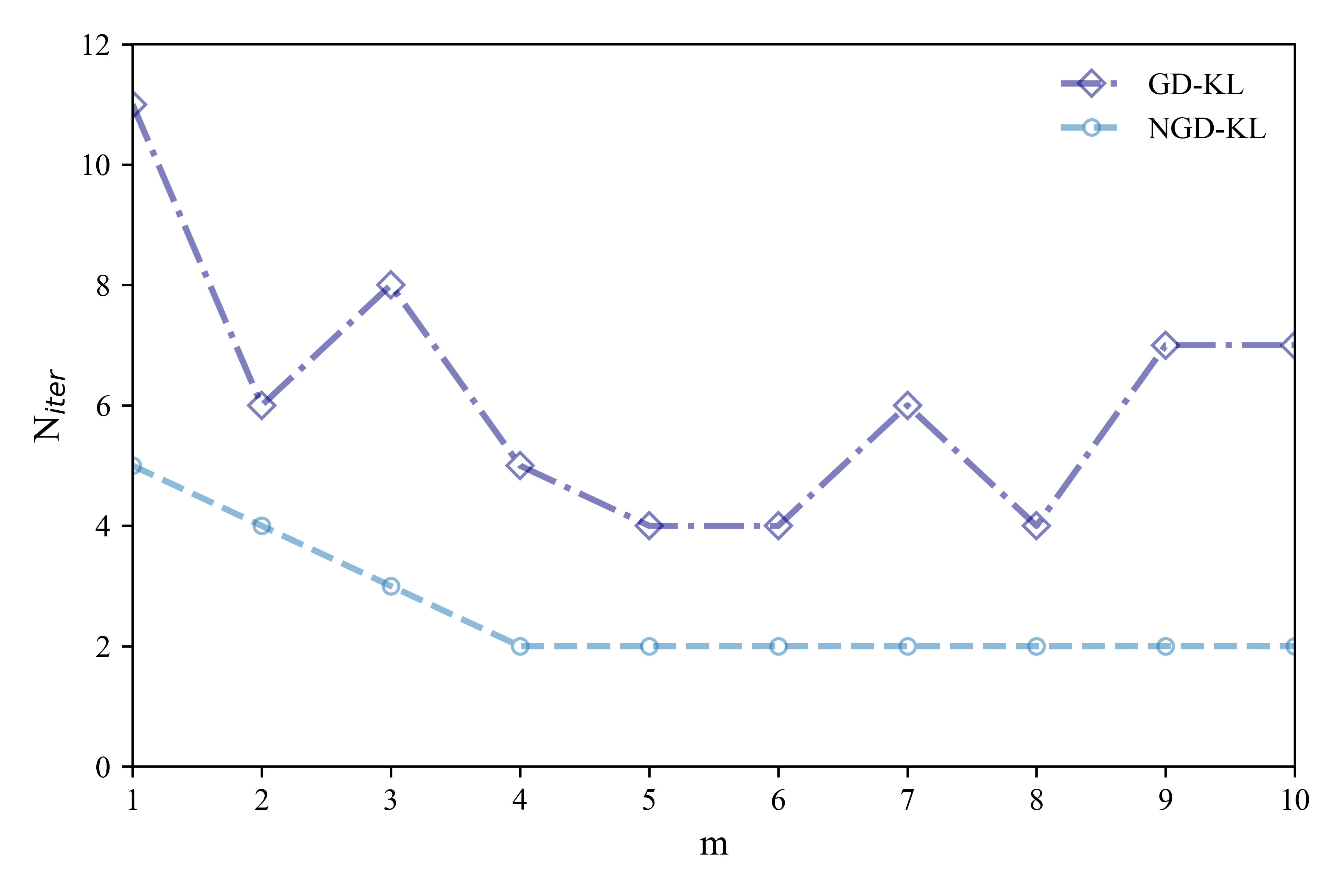}
\end{subfigure}
\begin{subfigure}{.49\textwidth}
    \includegraphics[width=\textwidth]{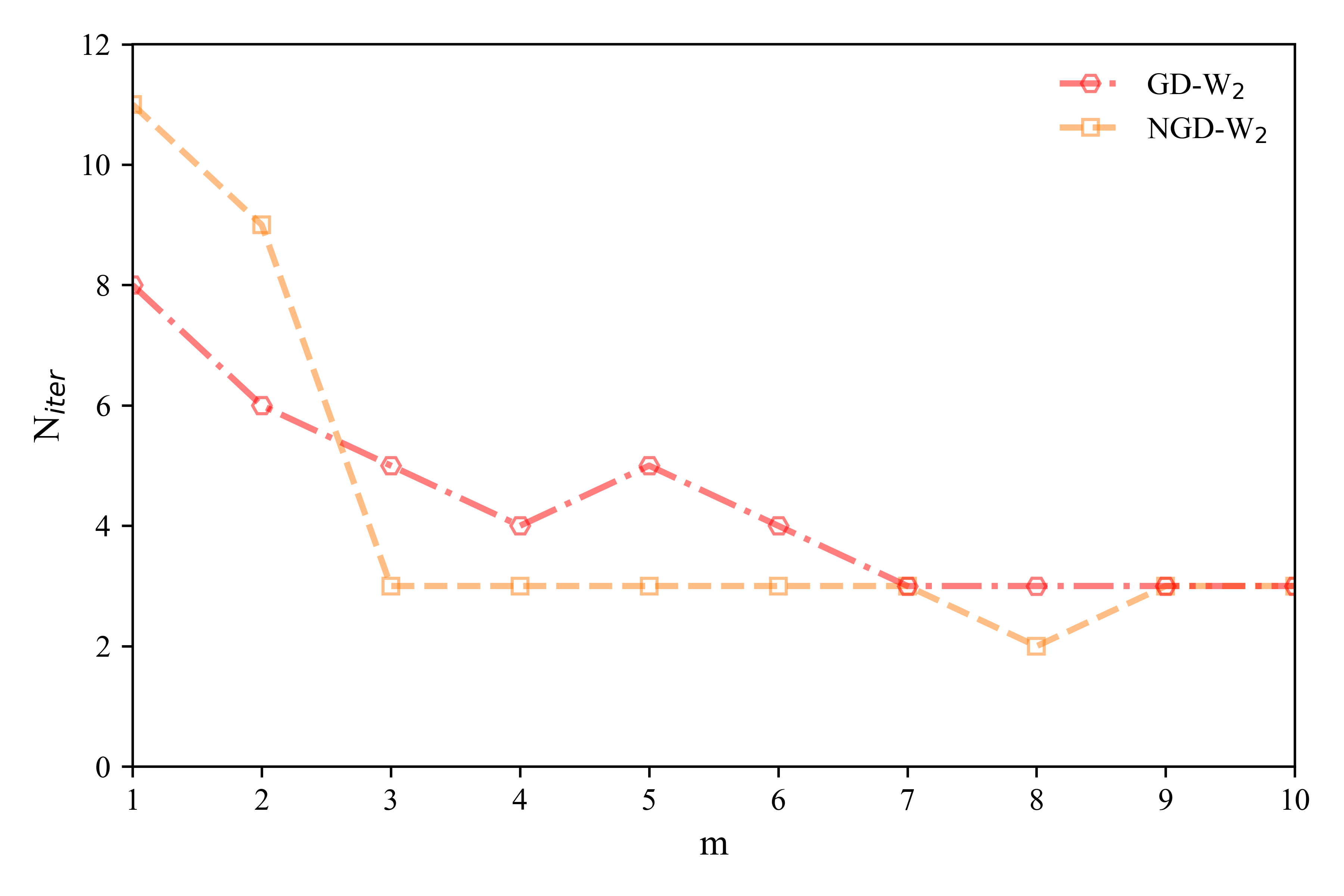}
\end{subfigure}
\end{center}
    \caption{Number of iterations to solve the minimization problem~\eqref{eq:mindd} at the $m$th assimilation step with GD (dash-dotted line) and NGD (dashed line) when the KL (left) and $W_2$ (right) discrepancies are used. The simulation parameter values are the same as in \cref{fig:pinitKL}.}
    \label{fig:NiterKL}
\end{figure}

\begin{figure}[htbp]
\begin{center}
\begin{subfigure}{.49\textwidth}
    \includegraphics[width=\textwidth]{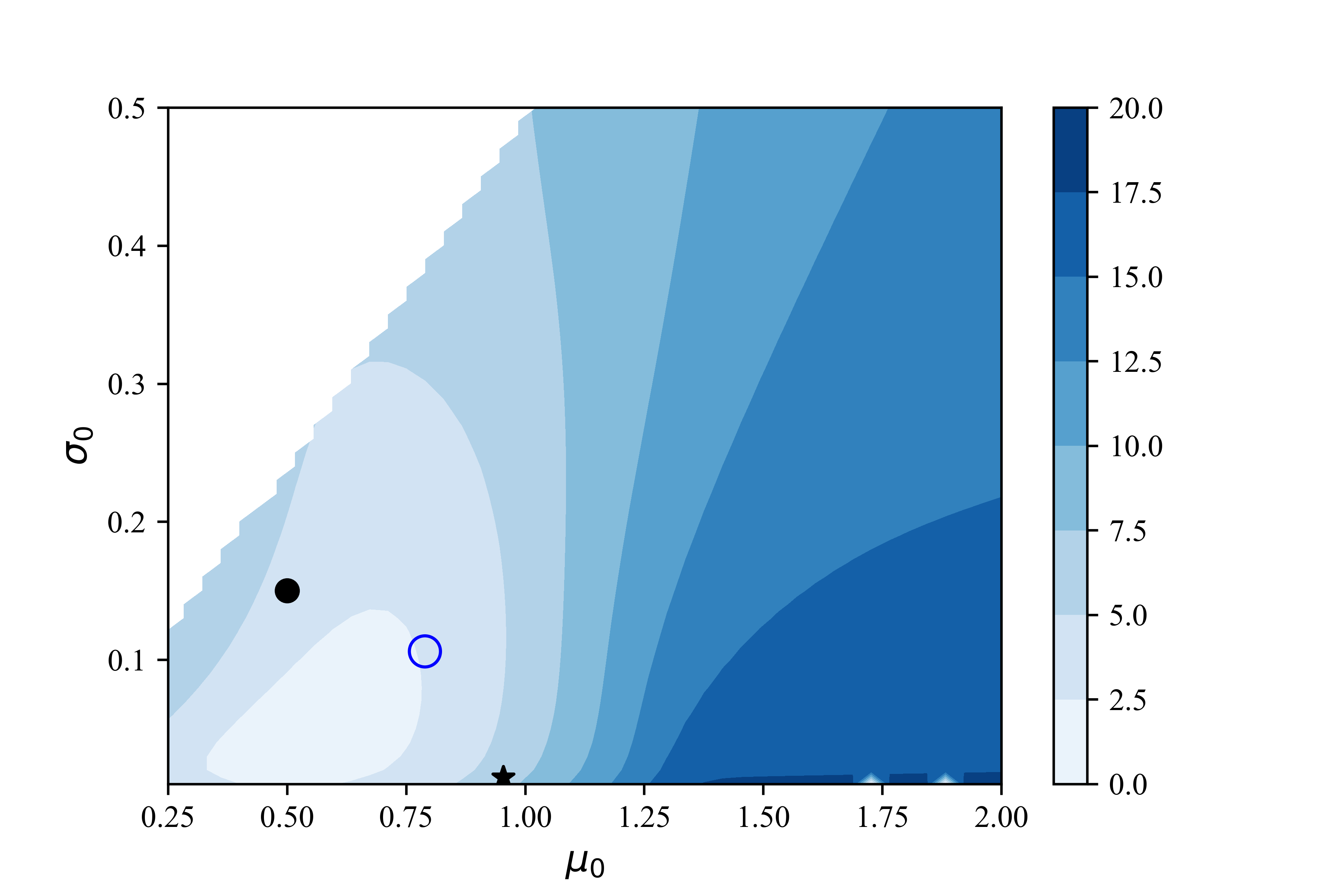}
\end{subfigure}
\begin{subfigure}{.49\textwidth}
    \includegraphics[width=\textwidth]{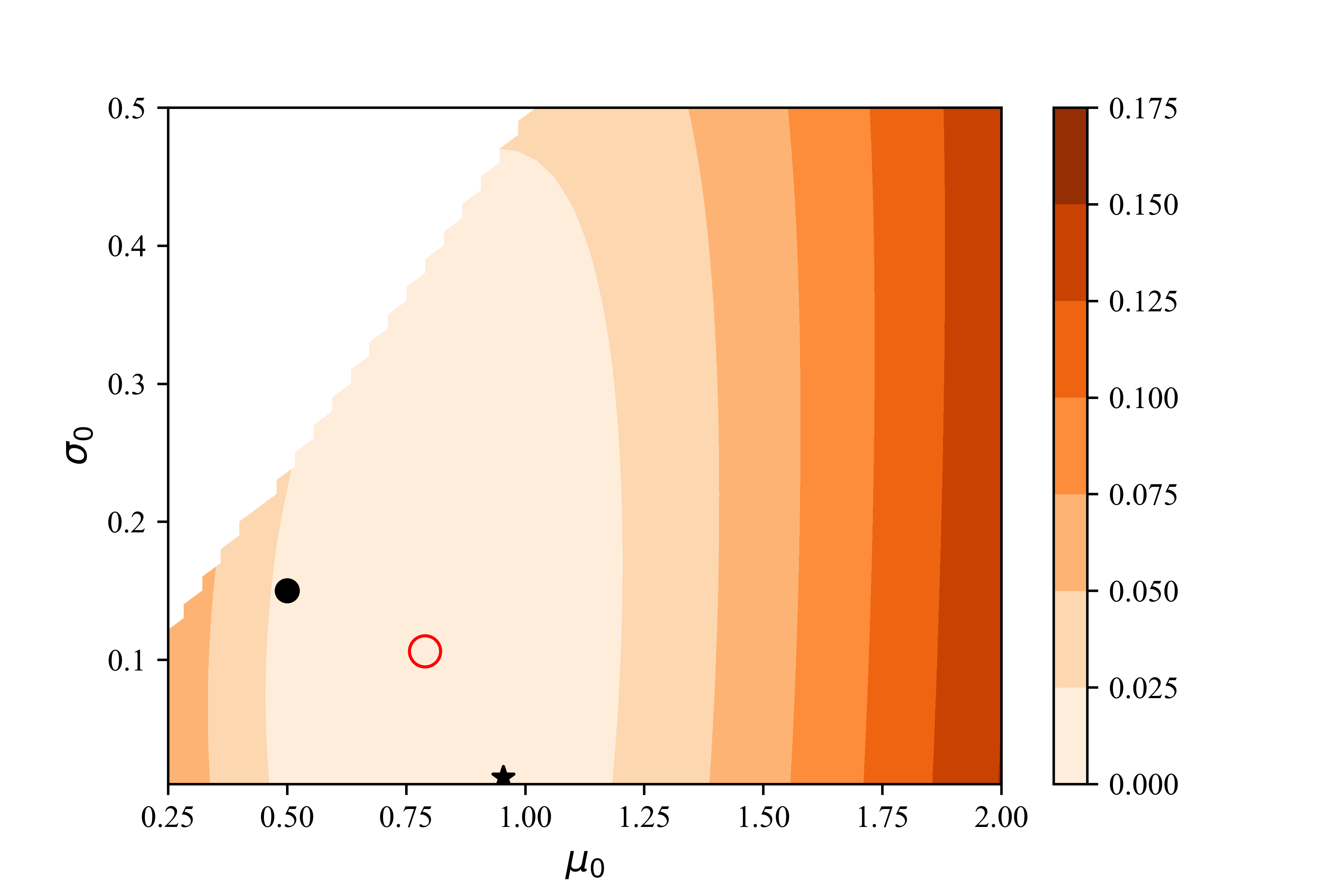}
\end{subfigure}
\end{center}
    \caption{Loss functions for KL (left) and W$_2$ divergence (right) for the random initial condition test case in the parameter space $\boldsymbol \varphi = \boldsymbol \varphi_0 = \{ \mu_0, \sigma_0 \}$ at $m=1$ (assimilation of the first measurement). The small (full) circle in both panels indicates the prior parameters $\boldsymbol \varphi^{(0)}$, the star indicates the true value $\{x_0^\star ,0 \}$. The large (empty) circle represents the first update of the parameters, $\boldsymbol \varphi^{(1)}$, obtained via NGD optimization. The blank region is to guarantee that $x_0$ is positive almost surely. The simulation parameter values are the same as in \cref{fig:pinitKL}. }
    \label{fig:losstc1}
\end{figure}

Figure~\ref{fig:NiterKL} shows the number of iterations, $N_\text{iter}$, it takes each of the  four DA methods to converge at each assimilation step of DA-MD. %Like before, the left panel refers to the KL loss whereas the right panel refers to the W$_2$ loss. 
If the KL divergence is used as a discrepancy metric, NGD converges in consistently fewer iterations than GD does; if the $W_2$ distance is used instead, then NGD and GD require on average the same number of iterations to converge. That is possibly because the $W_2$-induced loss function $\mathcal C (\boldsymbol \varphi) \equiv   W_2(f(X;\boldsymbol \varphi),\hat f(X))^2 / 2$ is smoother than its KL-induced counterpart (as shown in \cref{fig:losstc1} for $m=1$) and, hence, the availability of the analytical gradient $\nabla_{\boldsymbol\varphi} D$ is as helpful as the preconditioning. NGD is expected to become more beneficial when the loss functions is more sensitive to some parameters than to others, or when parameters vary in widely different ranges.

The total computational cost depends not only on the number of iterations, but also on the time required to compute the necessary terms at each iteration. Each NGD iteration is more expensive than GD's because it requires the calculation of a preconditioning matrix (i.e., the metric tensors for the geometry of the manifold), with a number of operations $\mathcal O(N_\text{par}^3)$ (see \cref{rem:cost}). When the evaluation of either the loss function or its gradient is computationally expensive, the computational time for GD would significantly increase.

As a final note, we found the calculation of the KL loss to be more sensitive to the initial guess, with poor choices of the prior often resulting in poor convergence of the DA-MD procedure. That is caused by the high sensitivity of the KL divergence to the shape of the distributions, especially when they are very far apart and/or very sharp, resulting in poor numerical accuracy of the integration \eqref{eq:KL}.

\subsection{Langevin equation with white noise}
\label{sec:problem2}

The dynamics of state variable $x(t)$ is described by a Langevin equation,
\begin{align}\label{eq:ODE2}
    \frac{\text d x}{\text d t} = - a(t) x(t), \qquad  x(0) = x_0^\star, 
\end{align}
where the statistically homogeneous (stationary) random process $a(t) = \mu_a + \sigma_a w(t)$ has mean $\mu_a$ and standard deviation $\sigma_a$, with $w(t)$ denoting standard Gaussian white noise. The initial state $x_0^\star \in \mathbb R^+$ is deterministic. The process $a(t)$ is almost surely positive, which ensures that $x(t)$ has a compact support $\Omega \subset \mathbb R^+$ and, hence, the information geometry induced by the $W_2$ distance is rigorously defined.

The single-point CDF $F(X;t)$ of the state variable $x(t)$ in~\eqref{eq:ODE2} satisfies exactly a CDF equation (\cref{app:MD})
\begin{subequations}\label{eq:MD_TC2}
\begin{equation}
    \frac{\partial F}{\partial t} - \mu_a X \frac{\partial F}{\partial X}  = \frac{1}{2} \frac{\partial }{\partial X} \left( \sigma_a^2 X^2 \frac{ \partial F}{\partial X} \right),
\end{equation}
subject to initial and boundary conditions
\begin{align}
    F(x;0) = \mathcal H(X - x_0^\star), \qquad F(X_{\text{min}};t) = 0, \qquad F(X_{\text{max}};t) = 1.
\end{align} 
\end{subequations}
%It can be rewritten in the form \eqref{eq:CDFeq_closed} with $\mathcal U = - \mu_a X$ and $\mathcal D = \frac{1}{2} \sigma_a^2 X^2 $. 
In this example, CDF $F$ is parameterized by $\boldsymbol \varphi = \{ \mu_a, \sigma_a \}$, which, as before, we make explicit by writing $F(X;t,\boldsymbol \varphi)$. The values of $\boldsymbol \varphi$ are refined by assimilating observations $\hat{\mathbf x}$. 

A physics-informed DNN (\cref{fig:fcnn}) serves as a surrogate model that approximates the solution of the CDF equation~\eqref{eq:MD_TC2}. The training set consists of the finite-volumes solutions~\cite{wheeler2005fipy} of \eqref{eq:MD_TC2} at selected points $(X,t)$, computed for a number of  different combinations of meta-parameters $\boldsymbol \varphi$. The details of this and other computations are provided in the opening of \cref{sec:problem}. %Minimization of \eqref{eq:losssurrogatemodel} is performed using the L-BFGS-B method, with a random initialization of $\mathbf w$ and $\mathbf b$. The construction of the DNN is done using Tensorflow \cite{abadi2016tensorflow}, taking advantage of the embedded automatic differentiation capabilities which yield all necessary ingredients: $f = \frac{\partial F}{\partial x}$, $\frac{\partial f}{\partial \varphi_i}$ and $\frac{\partial F^{-1}(Y)}{\partial \varphi_i}$ in $\nabla_\varphi D$ for the $i$ components of $\boldsymbol \varphi$. 
In this experiment, DNN %\textbf{[Right?]} Yes.
function approximation is considered satisfactory upon reaching a value of the loss function of {$4 \cdot 10^{-4}$}. This %The good accuracy achieved by the surrogate model in this case 
high accuracy enables the deployment of the DNN surrogate for both the analysis and forecast steps, further accelerating the information-geometric optimization of \eqref{eq:mindd} with the \texttt{Scipy} conjugate gradient routine. 

\begin{remark}
For more complex problems, it might be advantageous to use the surrogate model only for the approximation of the gradients, while retaining the finite-volume solution of the CDF equation for prediction. Alternatively, it might be necessary to construct a surrogate model for the local CDF at each assimilation step $m$, hence constructing a surrogate model for the CDF solution at time $t_m$ thus reducing the dimensionality of the input for the DNN.
\end{remark}

\begin{figure}[htbp]
\begin{center}
\begin{subfigure}{.49\textwidth}
  % include first image
  \includegraphics[width=\linewidth]{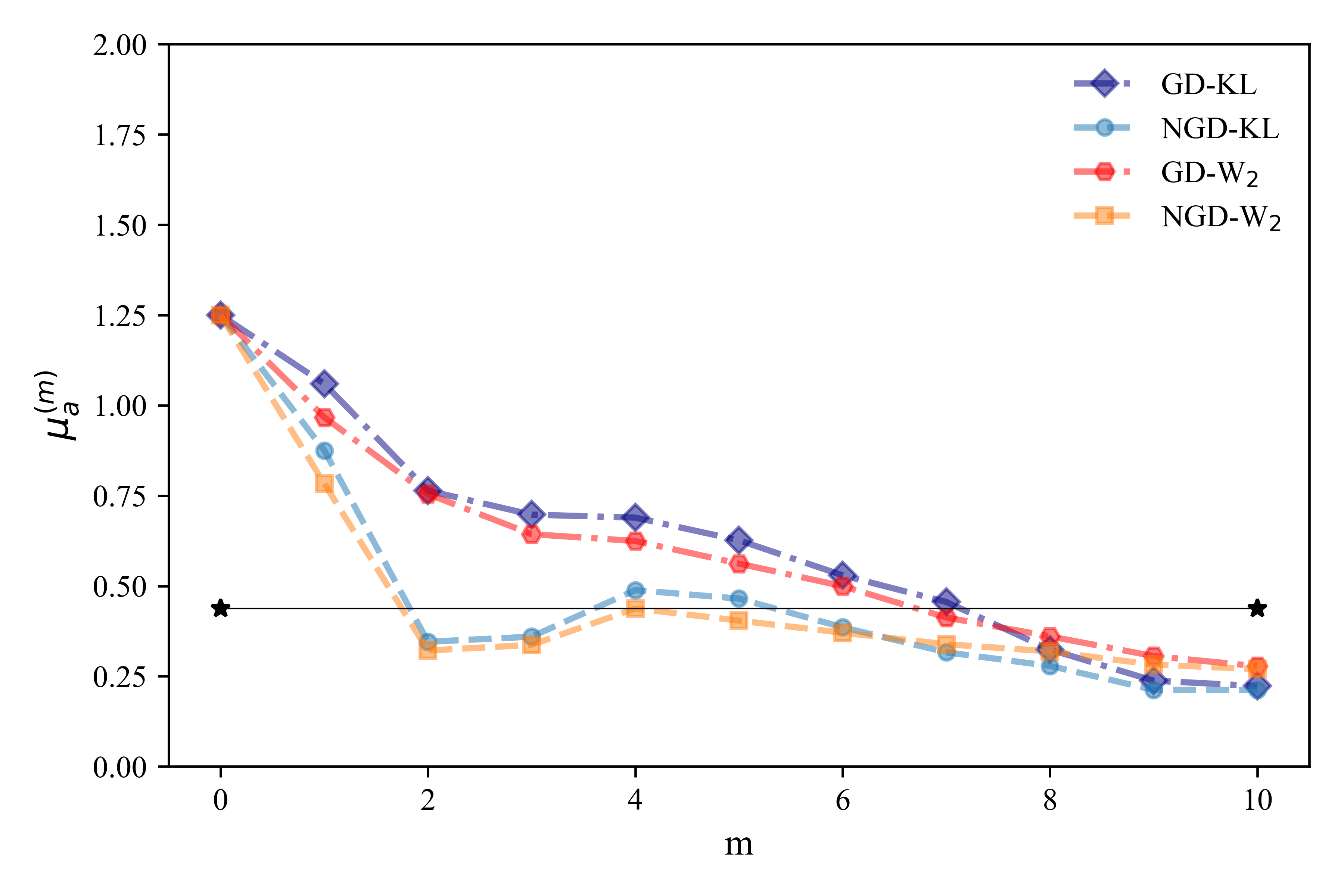}  
%   \caption{Put your sub-caption here}
%   \label{fig:sub-first}
\end{subfigure}
\begin{subfigure}{.49\textwidth}
  % include first image
  \includegraphics[width=\linewidth]{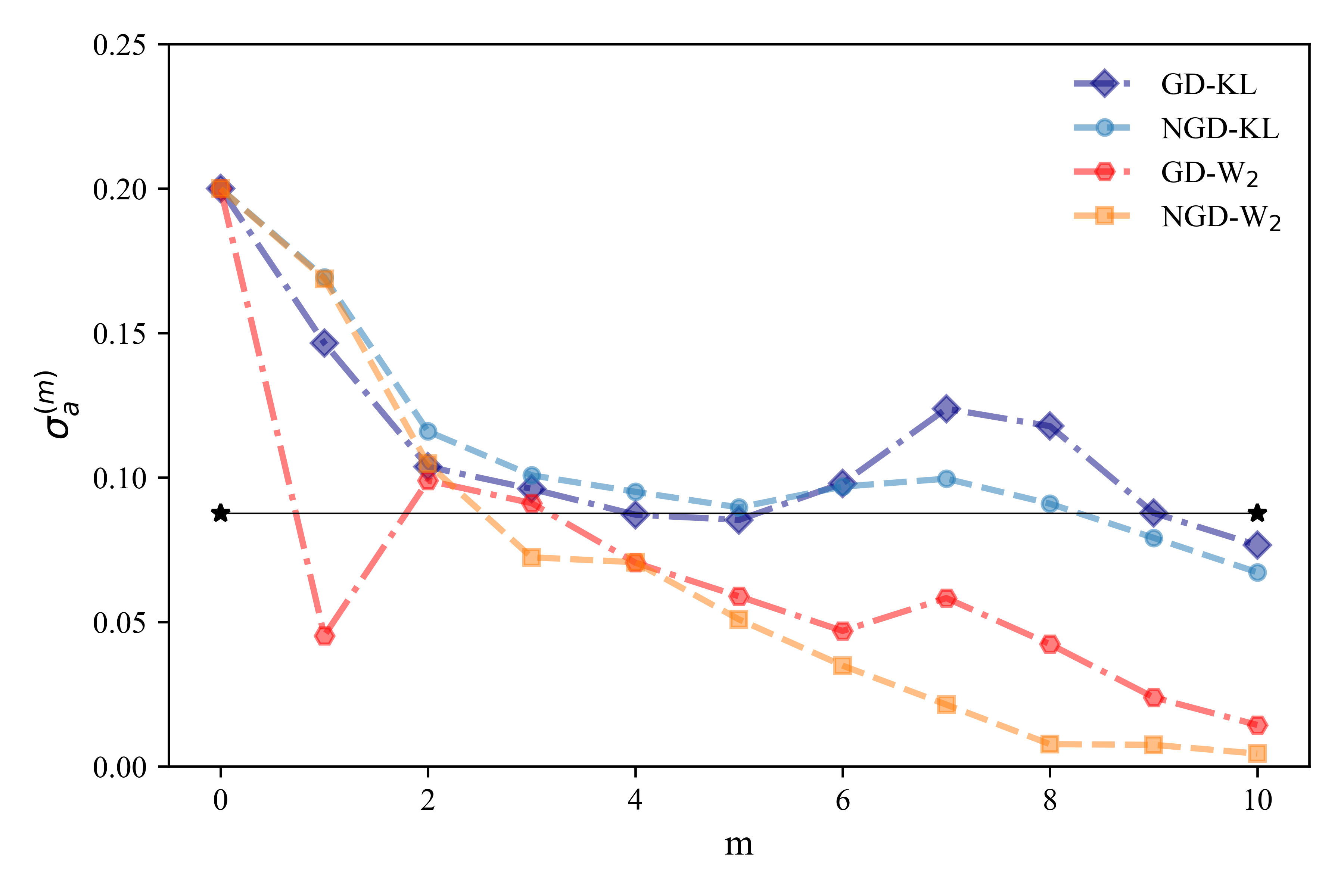}  
%   \caption{Put your sub-caption here}
%   \label{fig:sub-first}
\end{subfigure}
\end{center}
\caption{Estimation of meta-parameters $\boldsymbol \varphi = \{ \mu_a, \sigma_a\}$ for the Langevin equation with white noise. The parameters $\boldsymbol \varphi^{(m)}$ are plotted as function of the assimilation step $m$ for the four information-geometric optimization strategies: GD and NGD, for the KL and $W_2$ discrepancies. The simulation parameter values are set to $x_0^\star=1$, $\boldsymbol \varphi^\star = \{ 0.44, 0.088 \}$, $\boldsymbol \varphi^{(0)} = \{ 1.25, 0.2\}$ , $\sigma_\varepsilon = 0.1$, $N_\text{meas} = 10$, $t_{N_\text{meas}} = 2$, $\epsilon_{\text{KL}}=10^{-2}$, $N_T = 18850$, $N_B= 5632$, $N_I  = 1280$, and $N_R = 47872$. } 
\label{fig:metapars_case2}
\end{figure}

%Sequential Minimization of \eqref{eq:mindd} for both choices of the loss function targets the update of the meta-parameters $\boldsymbol \varphi = \{ \mu_a,\sigma_a\}$. 
Figure~\ref{fig:metapars_case2} shows the updated $\boldsymbol \varphi^{(m)}$ as function of the assimilation step $m$ for both KL and $W_2$ metrics of discrepancy, either taking advantage (NGD) or not taking advantage (GD) of the information-geometric structure of the statistical manifold of $F$. The starred values in this figure correspond to the statistical parameters used to generate the observations. All minimization algorithms converge to the exact mean $\mu_a$, whereas the identification of the standard deviation $\sigma_a$ is slightly more erratic. This is due to the inherent randomness of the physical process which calls for an improved data model, for example utilizing multiple observations at each observation time $t_m$.

\begin{figure}[htbp]
  \begin{center}
  \includegraphics[width=.5\linewidth]{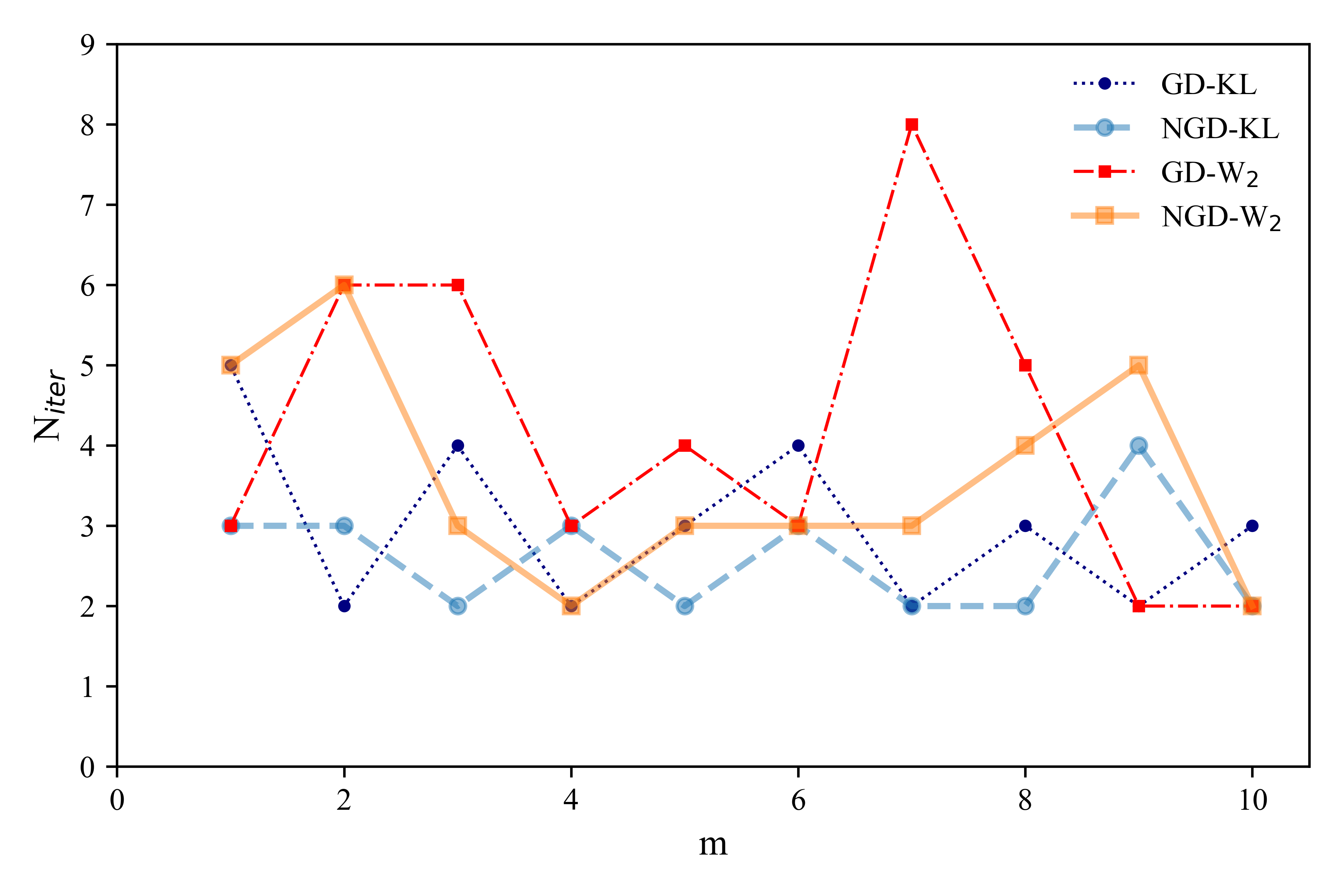}  
  \end{center}
\caption{Number of iterations per assimilation step $m$ for the four information-geometric optimization strategies: GD and NGD, for the KL and $W_2$ discrepancies. 
%A combination of line pattern and markers permits to discern among the different choices.  
The simulation parameter values are the same as in \cref{fig:metapars_case2}. }
\label{fig:Niters_case2}
\end{figure}

\begin{figure}[htbp]
\begin{center}
\begin{subfigure}{.49\textwidth}
  \includegraphics[width=\linewidth]{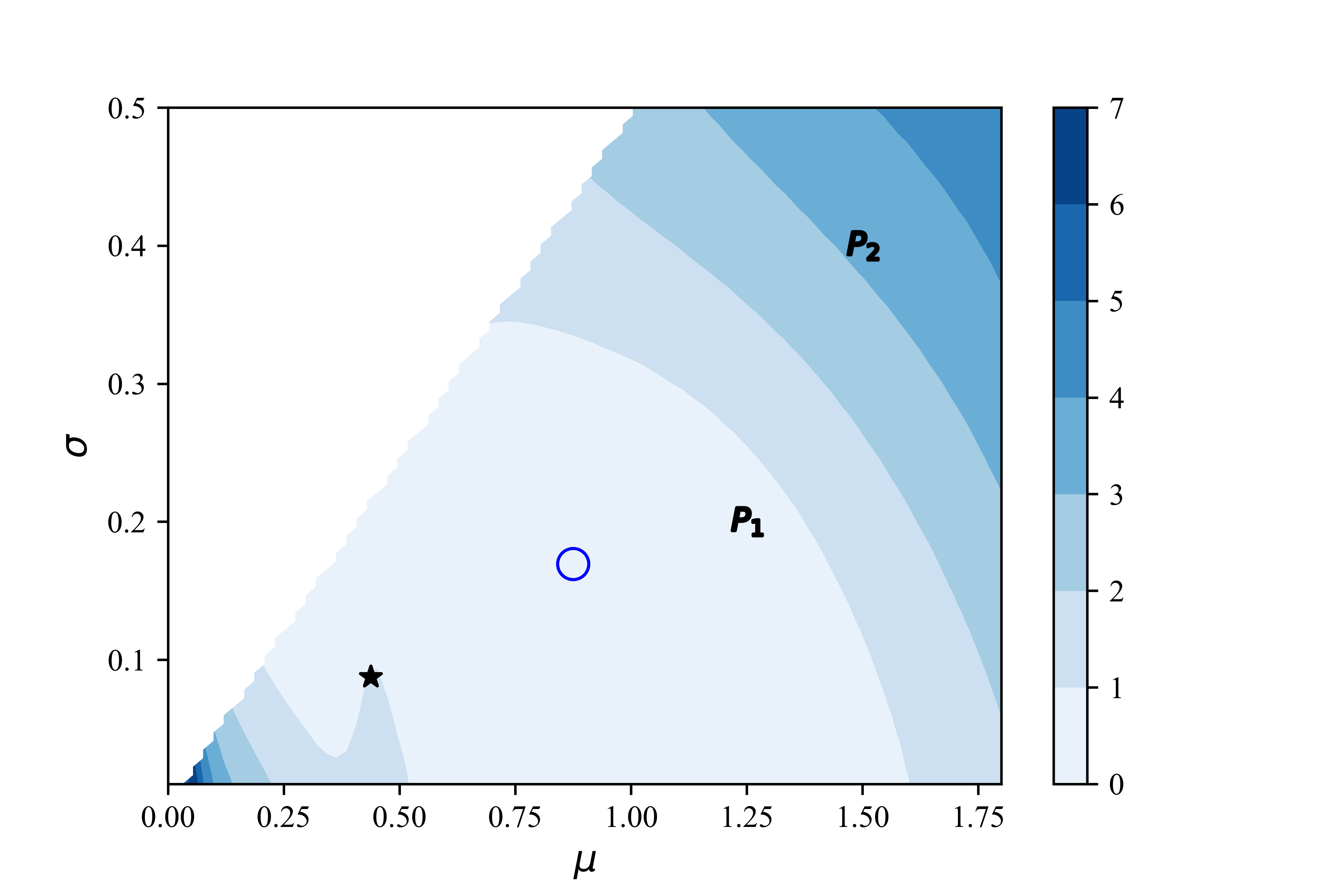}  
\end{subfigure}
\begin{subfigure}{.49\textwidth}
  \includegraphics[width=\linewidth]{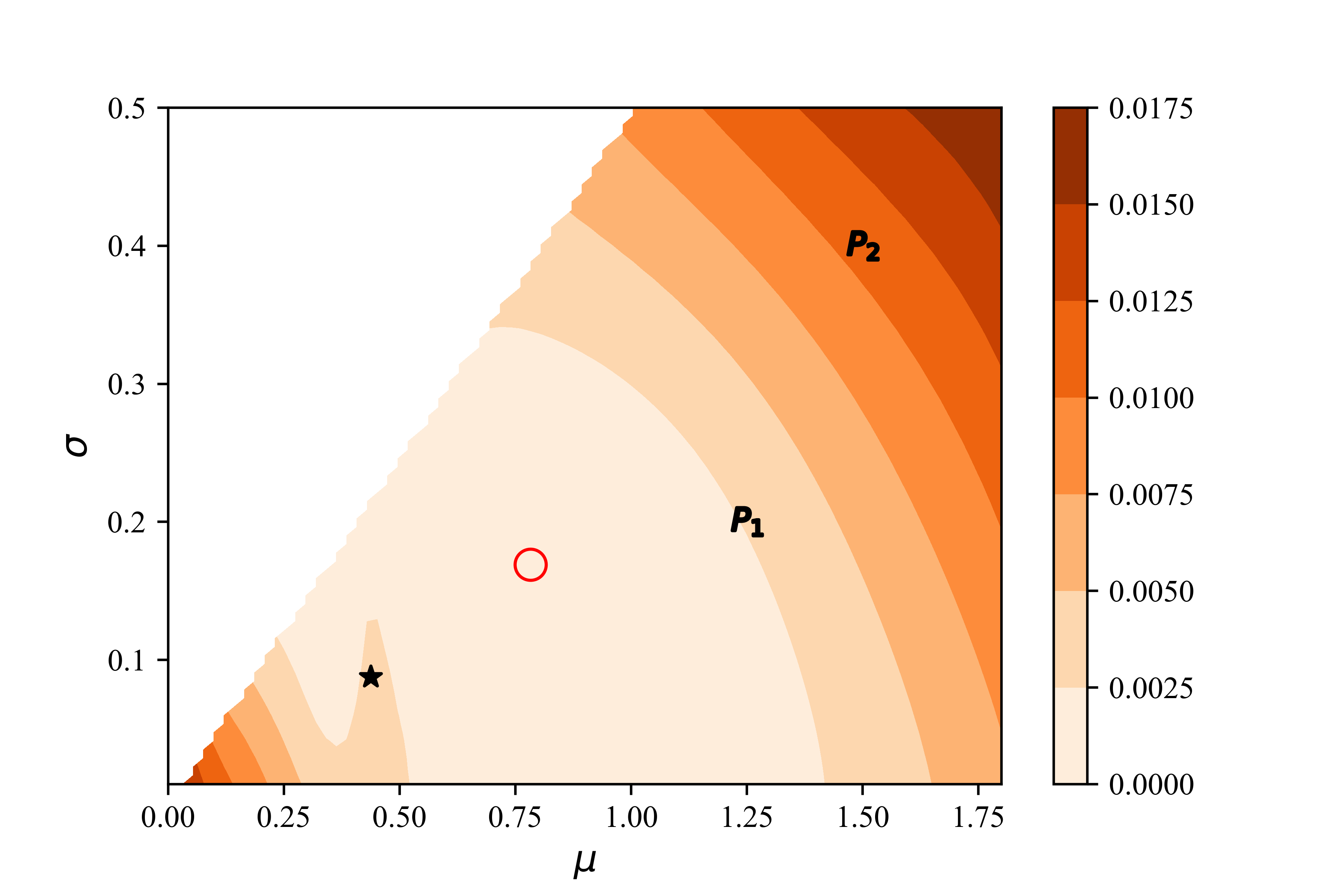}  
\end{subfigure}
\begin{subfigure}{.49\textwidth}
  \includegraphics[width=\linewidth]{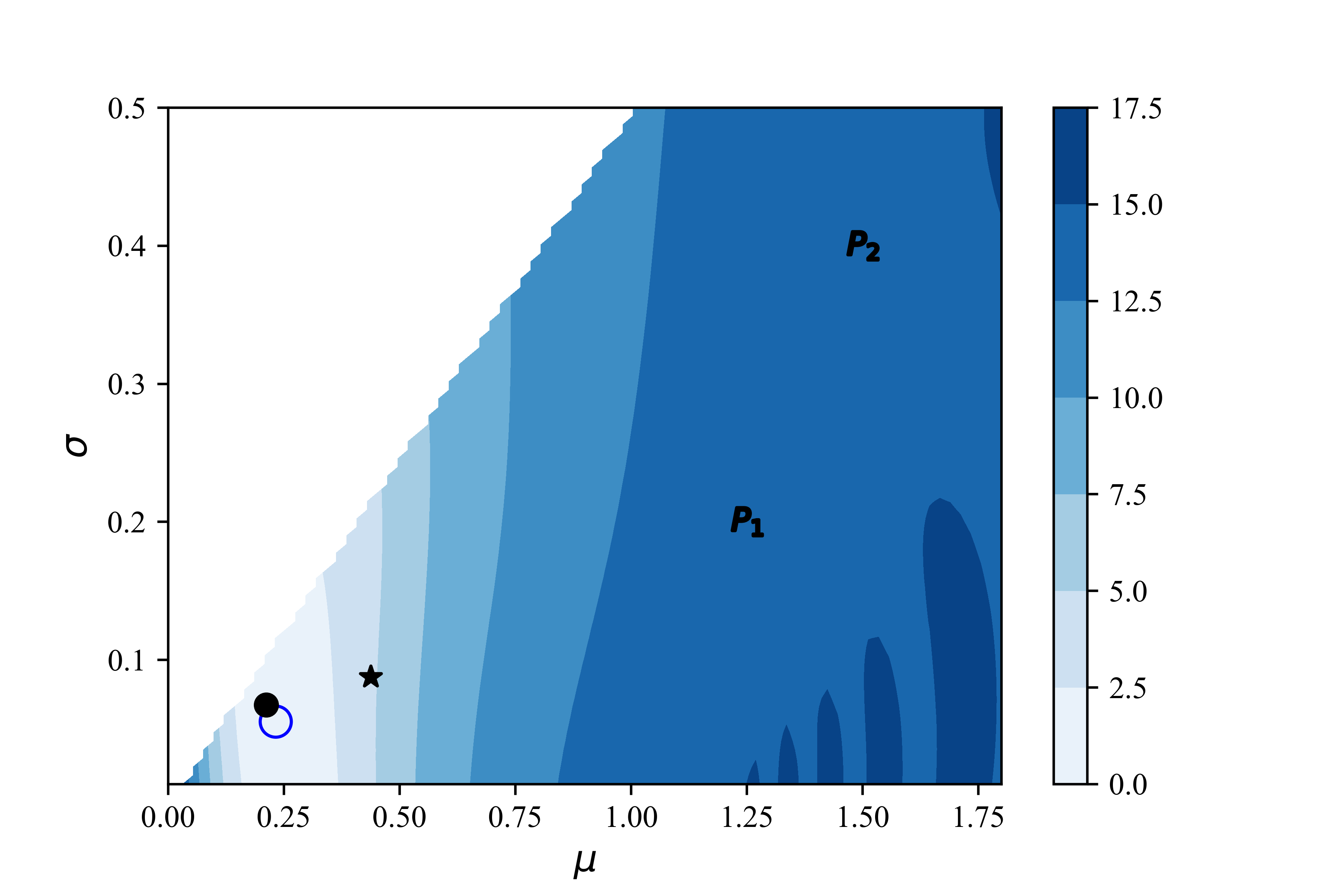}
\end{subfigure}
\begin{subfigure}{.49\textwidth}
  \includegraphics[width=\linewidth]{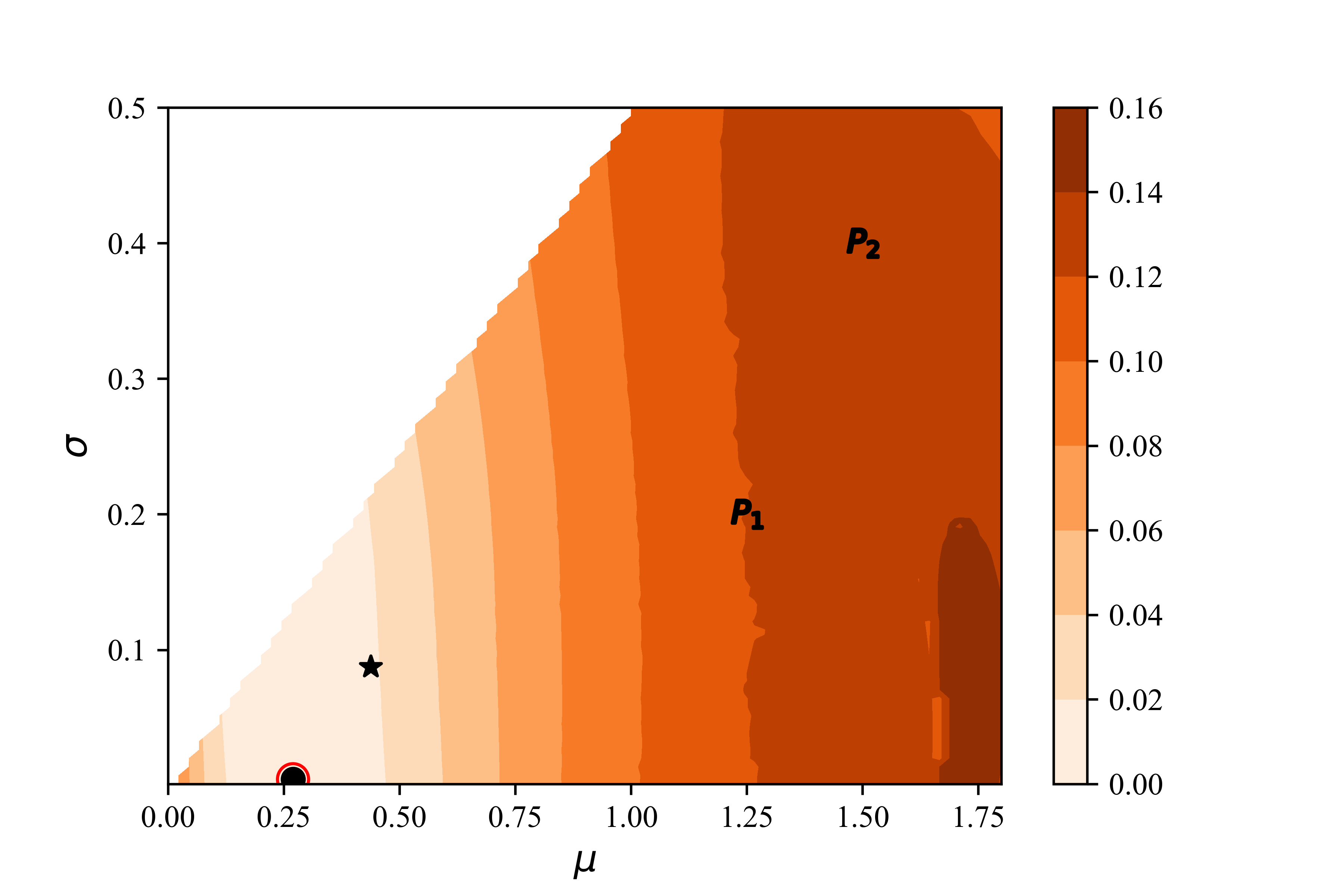}
\end{subfigure}
\end{center}
\caption{The KL (left column) and $W_2$ (right column) loss functions at the first ($m=1$, top row) and last ($m=N_\text{meas}$) steps of DA. The star indicates the true values of the meta-parameters (used to generate the synthetic reality). The points $P_1$ and $P_2$ indicates the priors $\boldsymbol \varphi^{(0)}$ for which the optimization of the KL loss function converges and fails to converge, respectively. The larger (empty) circles indicate the posterior parameters at the $m$th assimilation step, $\boldsymbol \varphi^{(m+1)}$, and the smaller (full) circles in the bottom row indicate $\boldsymbol \varphi^{(N_\text{meas}-1)}$. The blank region in all panels is to enforce almost surely the non-negativity of $a(t)$. 
The simulation parameter values are the same as in \cref{fig:metapars_case2}.}
\label{fig:losses_case2}
\end{figure}

Like before, the number of iterations over the assimilation time window is smaller for NGD than for GD for both choices of loss function (\cref{fig:Niters_case2}), albeit the difference is not as pronounced. The physics-driven parameterization of the statistical manifold yields an isotropic geometry of the loss function in the search area, which reduces the benefits of preconditioning.
This is shown in \cref{fig:losses_case2}, where the KL and $W_2$ loss functions are plotted, at the first and last assimilation steps, as function of the meta-parameters $\boldsymbol \varphi$, also highlighting the true solution and the prior location.\footnote{The loss functions at assimilation step $m=1$ are obtained using the initial $\boldsymbol \varphi^{(0)}$ for the calculation of the observational PDF/CDF, whereas the loss functions at  assimilation step $N_\text{meas} = 10$ are computed using $\boldsymbol \varphi^{(N_\text{meas}-1)}$ for the prior obtained using either NGD-KL or NKD-W$_2$. The initial guess of the prior $\boldsymbol \varphi^{(0)}$ is the same for both KL and $W_2$ metrics (P$_1$ in the Figure), and yields similar outcomes in terms of identification of the meta-parameters, as illustrated above.} Although superficially similar throughout the assimilation process, the minor differences in the topology of the KL and W$_2$ loss functions are enough to prevent convergence for the KL loss function for a slightly worst choice of the prior (P$_2$ in the Figure), which results in divergence of the DA-MD procedure for both GD and NGD. This is because the KL divergence is more sensitive to numerical errors in the calculation of the integrals, especially for sharp or non-overlapping distributions, which mislead the direction of the search.

\begin{figure}[htbp]
\begin{subfigure}{.5\textwidth}
  \centering
  % include first image
  \includegraphics[width=\linewidth]{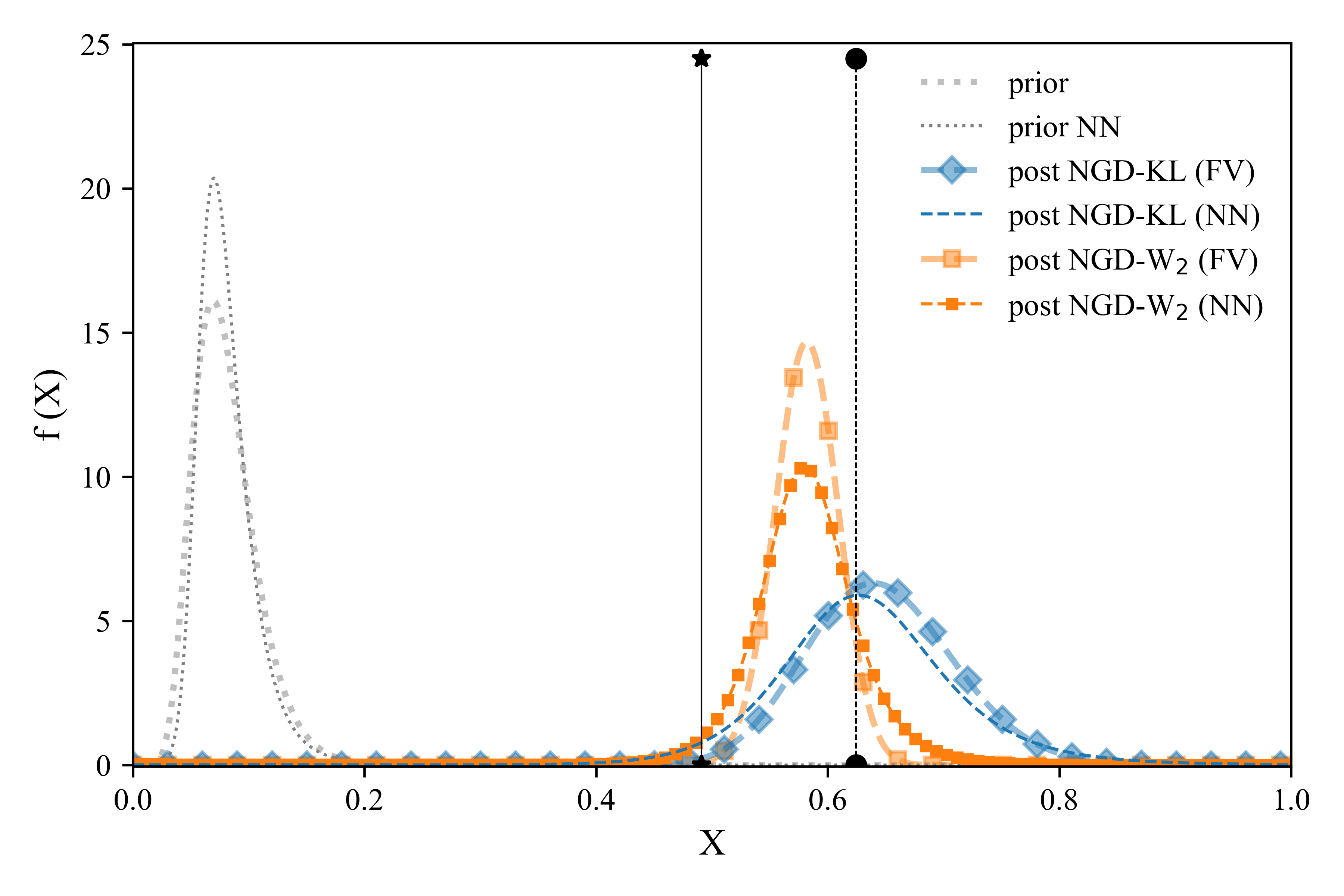}  
%   \caption{Put your sub-caption here}
%   \label{fig:sub-first}
\end{subfigure}
\begin{subfigure}{.5\textwidth}
  \centering
  % include second image
  \includegraphics[width=\linewidth]{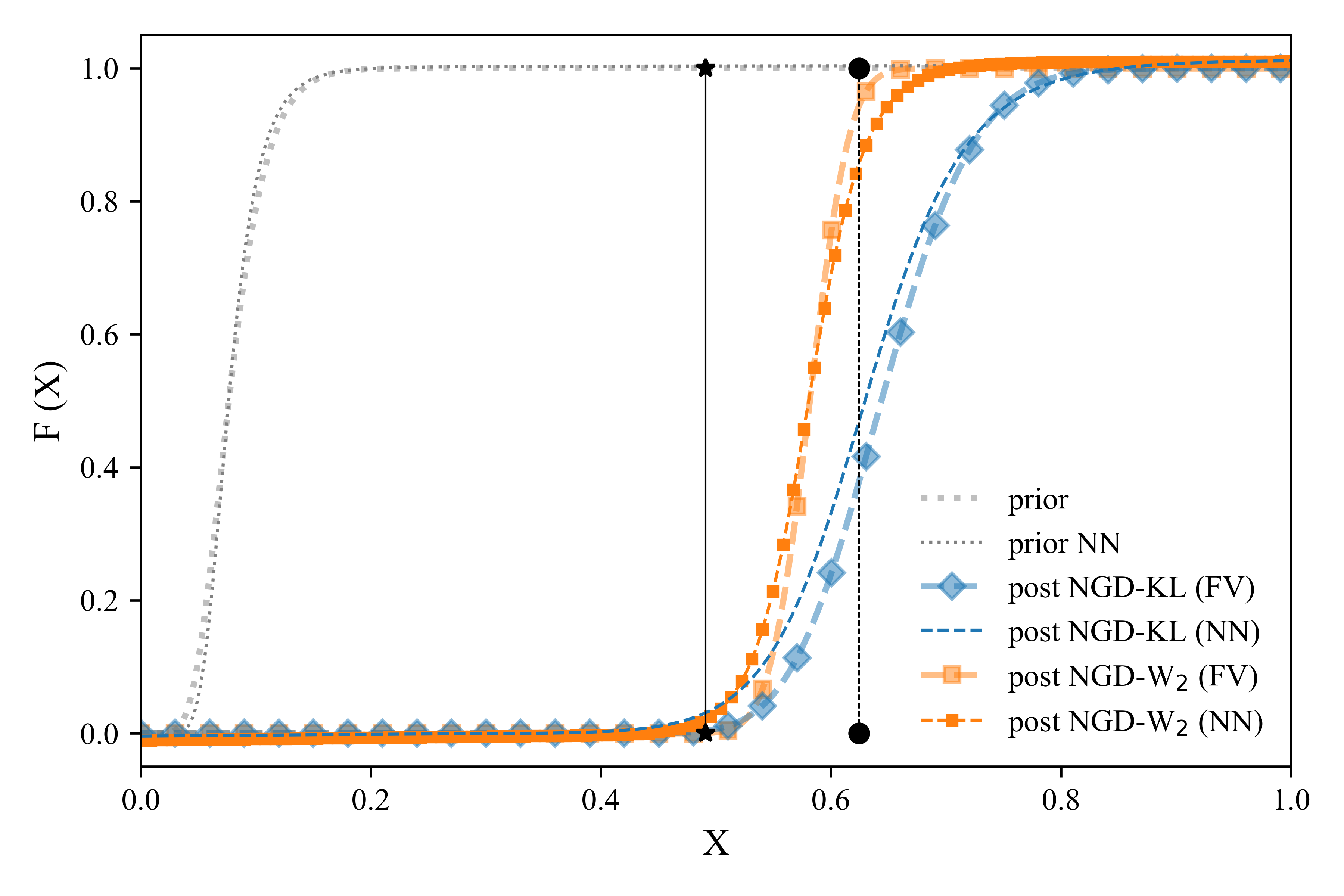}  
%   \caption{Put your sub-caption here}
%   \label{fig:sub-second}
\end{subfigure}
\caption{Prior and posterior distributions (PDFs on the left, and corresponding CDFs on the right) at time $t_{N_\text{meas}}$ obtained via NGD minimization for the KL and $W_2$ loss functions, with either the FV solution of the CDF equation~\eqref{eq:MD_TC2} or its NN surrogate. Black stars and circles mark the exact value $x(t_{N_\text{meas}})$ and its noisy observation $x_{N_\text{meas}}$, respectively. The simulation parameter values are the same as in \cref{fig:metapars_case2}.  }
\label{fig:posteriors_case2}
\end{figure}

The posterior NGD parameters $\boldsymbol \varphi^{(N_\text{meas})}$ are used to compute the posterior CDF and PDF of $x(t)$ in \cref{fig:posteriors_case2}.  NGD yields accurate posteriors, with the $W_2$ optimization~\eqref{eq:WGD} performing better than the KL optimization~\eqref{eq:FGD}. In order to highlight the accuracy of the DNN surrogate model, we show the finite-volume solution of the CDF equation~\eqref{eq:MD_TC2} with $\boldsymbol \varphi = \boldsymbol \varphi^{(N_\text{meas})}$ and its corresponding PDF computed via numerical differentiation, and their DNN-based counterparts. In agreement within~\cite{chen2018wasserstein}, we found the $W_2$ minimization to be more robust to the choice of the prior.

\begin{remark}
An additional advantage of the $W_2$ loss function stems from its reliance on a CDF rather than a PDF that enters the KL loss function. CDFs are smoother and easier to compute as a solution of the CDF equation than PDFs, which are obtained by solving the PDF equation. This facilitates the generation of a training set and the training of a surrogate model. On the other hand, approximation of the solution to a CDF equation with a DNN surrogate possesses a potential challenge for the $W_2$ optimization, since~\eqref{eq:mindd} calls for invertible surrogate models. We overcome this difficulty by 
selecting a special structure for the DNN that guarantees automatic inversion, as detailed in \cref{sec:method}. 
\end{remark}

% The chosen setup of FCNN makes the calculation of W2 easier once a NN surrogate is available because it removes the necessity of an invertible neural network of the usual type (with just one output). overall W2 is preferable because of that, and because it is more robust. in our case it also yields with the same computational time better overall accuracy of the result. The dicrepancy in terms of computational time is destined to increase the more expensive it becomes to compute function evaluations (not always reliable to use the surrogate model for that), for example the experiments we did in the other paper go back to the solution of the PDE when recomputing the observational CDF at each iteration).
%, especially when the prior is far from the true value of the distribution. 

The computational cost of the different optimization strategies depends on the number of iterations (\cref{fig:Niters_case2}); on the computational cost per iteration; and, in case of information-geometric optimization, on the cost of computing the tensor metrics. Since the function- and gradient-evaluations for this example are not expensive, the computational gain of having a smaller number of evaluations is not significant, and it is compensated by the additional cost of the calculation of the preconditioning matrices. %The differential would emerge in more complex cases. 

\subsection{Langevin equation with colored noise}
\label{sec:problem3}

The dynamics of state variable $x(t)$ is described by~\eqref{eq:ODEphysical} with $s(x(t); w(t), \boldsymbol \theta ) \equiv - a(t) x(t)$, where
$a(t) = \mu_a + w(t)$ with $\mu_a \in \mathbb R^+$, and $w(t)$ is the derivative of an Ornstein–Uhlenbeck process characterized by the exponential auto-covariance function
\[
C_w(t,\tau) = \frac{\sigma_a^2}{2 \theta_a} \left[ \text{e}^{ - \theta_a | t - \tau | } + \text{e}^{ - \theta_a (t+\tau) } \right],
\]
with parameters $\sigma_a$ and $\theta_a \in \mathbb R^+$. By construction, the latter is also the auto-covariance function of $a(t)$, $C_w(t,\tau) = C_a(t,\tau)$. 
% \begin{comment}
% We are solving the stochastic ODE
% \begin{equation}
%     d x_t = - \mu_a x_t \text d t - x_t d OU_t
% \end{equation}
% where $OU_t$ is the Ornstein-Uhlenbeck process, which is characterized by
% \begin{equation}
%     d OU_t = - \theta OU_t \text d t + \sigma d W_t,
% \end{equation}
% with positive parameters $\sigma$ and $\theta$, and where $W_t$ is the standard Wiener process. Using JITCSE, we define $x_t = y^{(0)}$, and $\text d OU_t / dt = y^{(1)}$, so that the stochastic system to be solved reads
% \begin{align}
%     & \frac{\text d y^{(0)}}{dt} = - mu_a y^{(0)} - y^{(0)} y^{(1)} \\
%     & \frac{\text d y^{(1)}}{dt} = - \theta_a y^{(1)} + \sigma \eta
% \end{align}
% where $\eta$ is the derivative of the standard Wiener process. It can be shown that the mean and the covariance of the OU process are as written in the main text https://en.wikipedia.org/wiki/Ornstein–Uhlenbeck_process#In_physical_sciences
% \end{comment}
%; and the deterministic but unknown constants $\mu_a, \sigma_a, \theta_a \in \mathbb R^+$ are the mean, standard deviation, and the reciprocal of the correlation time of $a(t)$, respectively. 
Taking the initial state $x_0$ to be deterministic, the stochastic solution of this problem depends on three meta-parameters $\boldsymbol \varphi = \{ \mu_a, \sigma_a, \theta_a\}$. One realization of this solution, drawn from the distribution with the ``true'' meta-parameters $\boldsymbol \varphi^\star$, serves as ground truth for which observations $\hat{\mathbf x}$ are constructed in accordance with~\eqref{eq:data}.

%--- Moreover, by choosing $x_0$ with compact support $\Omega_0 \subset \mathbb R^+$ in the first scenario, and positive in the second and third scenario, and positive $a(t)$, we ensure that $x(t)$ has a compact support $\Omega \subset \mathbb R^+$. {This is to make sure that the information geometry induced by the W$_2$ divergence is rigorously defined}. ---

%We start by recasting this problem in terms of the standard Gaussian white noise $\xi(t)$ as
%\begin{equation}
%\begin{aligned}
%& \frac{\text d x}{\text d t} = - \mu_a x(t) - w(t) x(t)  \\
%& \frac{\text d w}{\text d t} = - \theta_a w(t) - \sigma_a \xi(t).
%\end{aligned}
%\end{equation}

We show in \cref{app:problem3} that the CDF $F(X;t)$ of $x(t)$ satisfies the CDF equation~\eqref{eq:CDFeq_closed} with 
\begin{align}\label{eq:UDtest3}
    \mathcal U(X,t;\boldsymbol\varphi) = -\mu_a X + X \! \int_0^t \!\! C_w(t,\tau) \text d \tau \quad \text{and} \quad
    \mathcal D(X,t;\boldsymbol\varphi) = X^2 \! \int_0^t \!\! C_w(t,\tau) \text d \tau.
\end{align}
The FV solution of this equation and its DNN surrogate are used to assimilate observations $\hat{\mathbf x}$ via our information-geometric DA-MD framework. Similar to the case of white noise (\cref{sec:problem2}), we found the KL-based implementation of DA-MD to be less robust to the choice of the prior. Hence, only the $W_2$-based results are displayed below.

Figure~\ref{fig:metapars_case3} exhibits the convergence of the meta-parameters $\boldsymbol\varphi$ as function of the data assimilation step $m$. Since the $W_2$ loss function is relatively insensitive to the third meta-parameter $\theta_a$, we present the convergence results for $\tilde \sigma = \sqrt{ \sigma^2 / (2 \theta_a)}$ instead.\footnote{This lack of sensitivity reflects the challenge of inferring the correlation length, $1/\theta_a$, from observations over a time window spanning only two true correlation lengths, $1/\theta_a^\star$.} Both GD-W$_2$ and NGD-W$_2$ converge after assimilation of about 20 observations, which are generated every $\Delta t = 0.055$. NGD converges, for the given combination of observations and the prior, in fewer iterations over the assimilation window (\cref{fig:metapars_case3}d) than GD. 

\begin{figure}[htbp]
\begin{center}
\begin{subfigure}{.49\textwidth}
  \includegraphics[width=\linewidth]{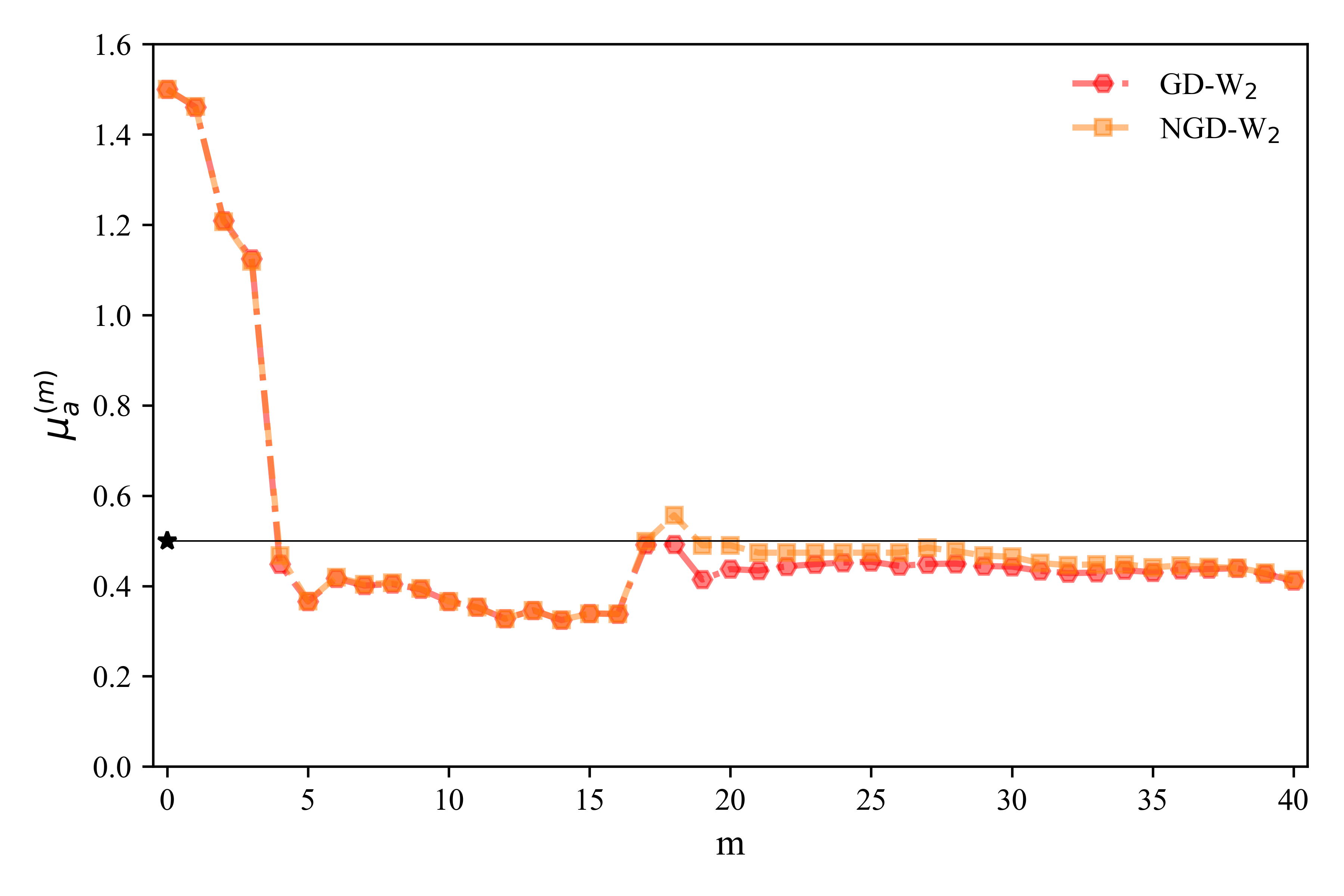}  
\end{subfigure}
\begin{subfigure}{.49\textwidth}
  \includegraphics[width=\linewidth]{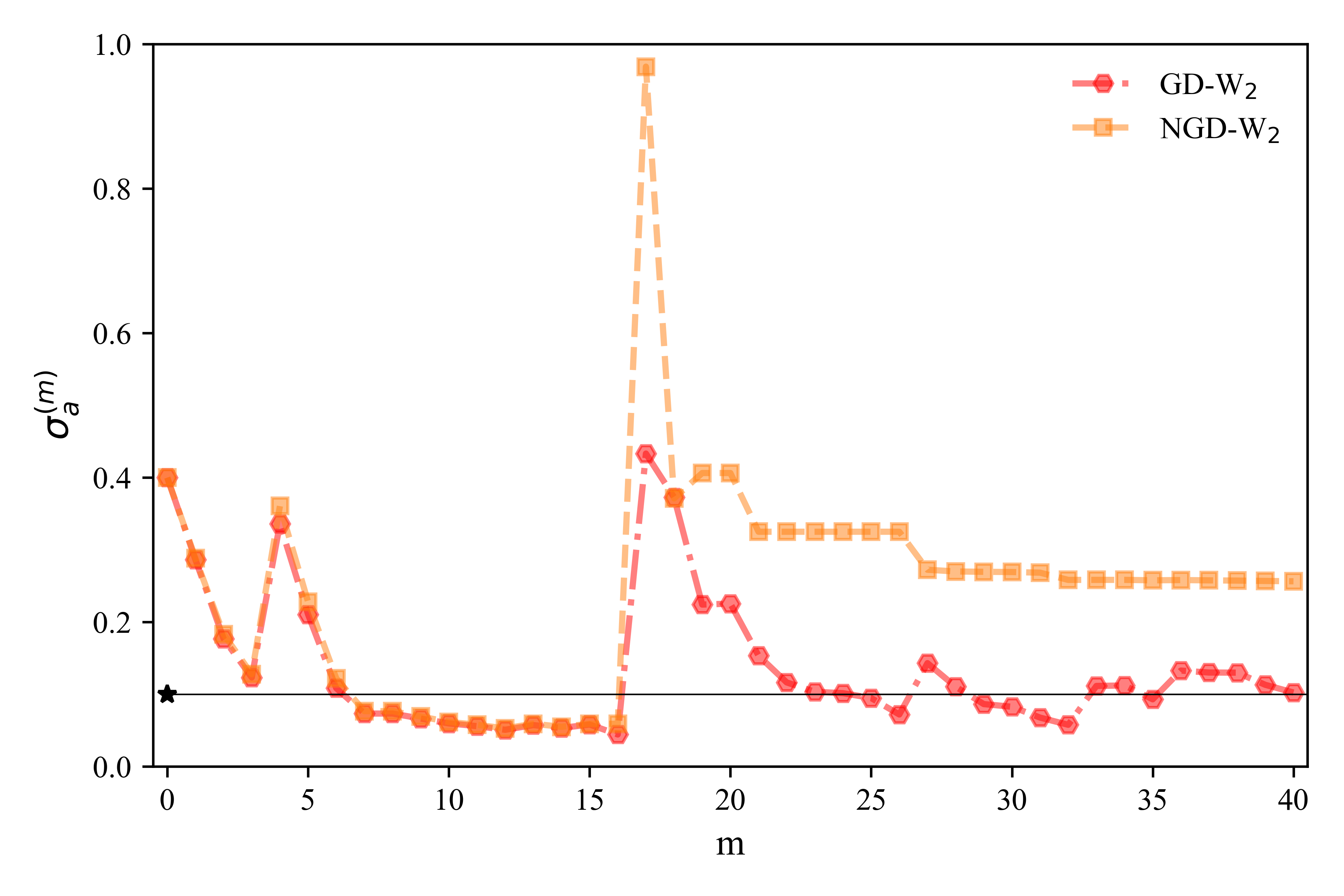}
\end{subfigure}
\begin{subfigure}{.49\textwidth}
  \includegraphics[width=\linewidth]{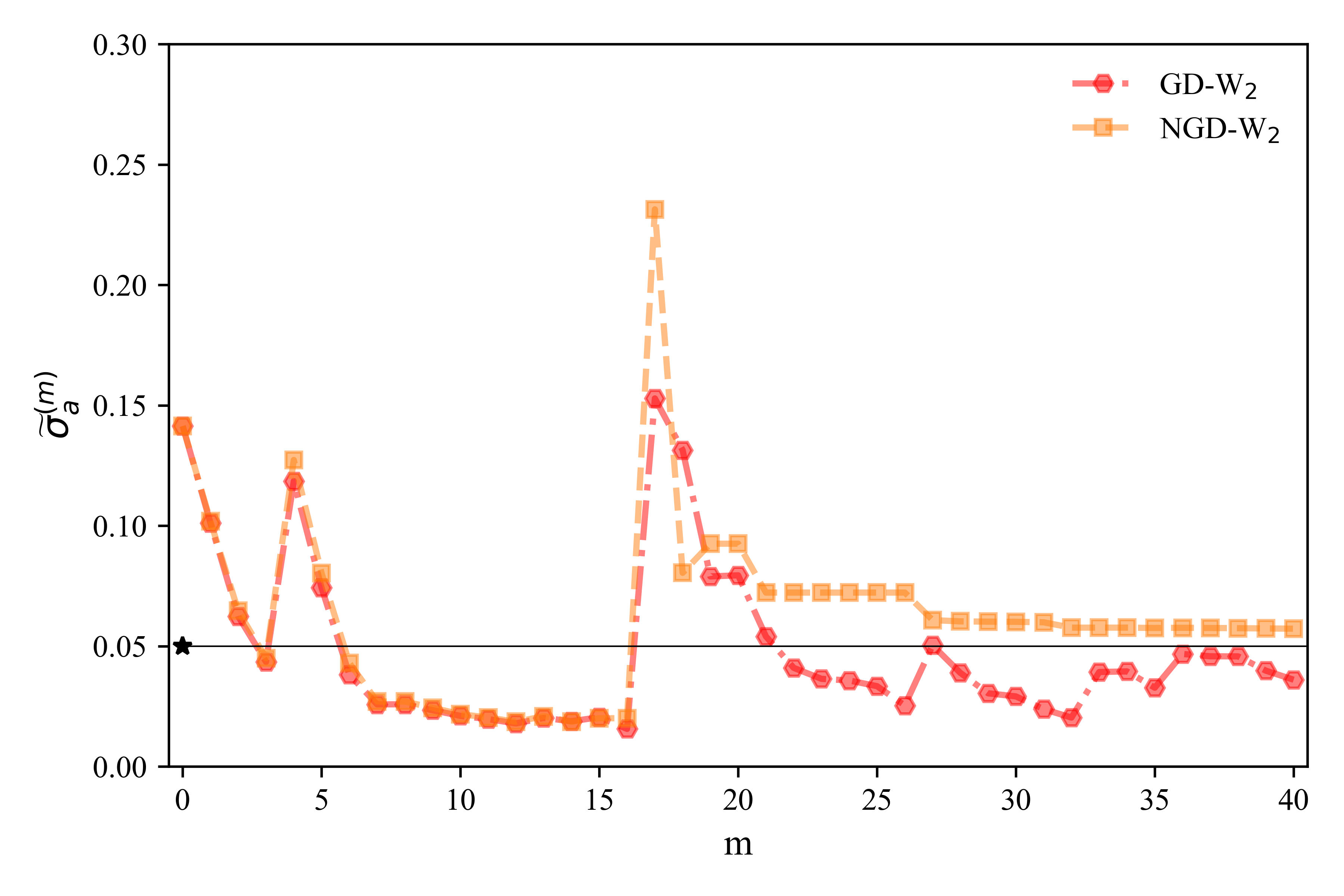}
\end{subfigure}
\begin{subfigure}{.49\textwidth}
  \includegraphics[width=\linewidth]{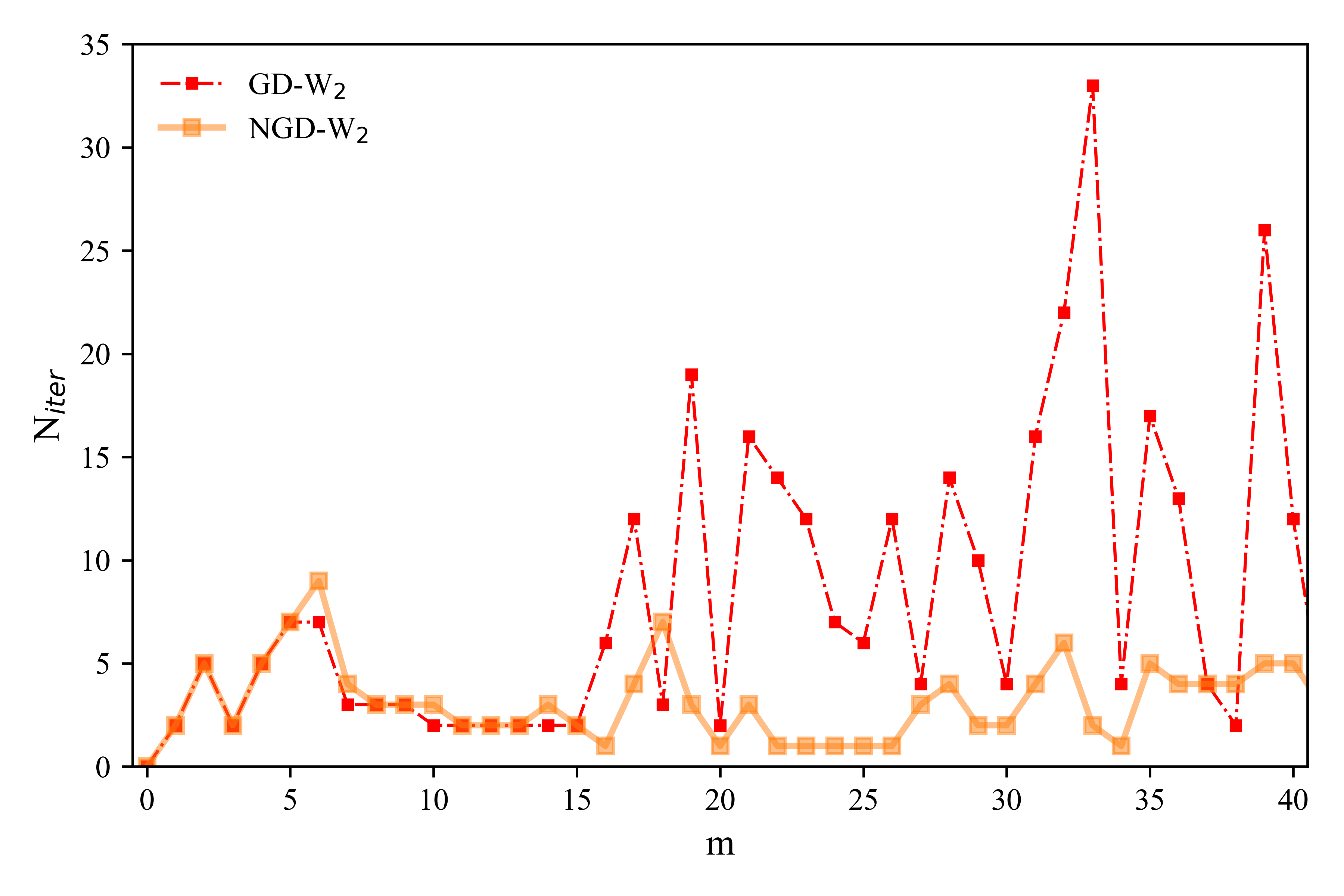}  
\end{subfigure}
\end{center}
\caption{Estimation of meta-parameters $\boldsymbol \varphi = \{ \mu_a, \sigma_a, \tilde \sigma = \sqrt{ \sigma_a^2/(2\theta_a) } \}$, as  function of the assimilation step $m$, with GD and NGD for the $W_2$ loss functions. The bottom right panel shows the number of iterations per assimilation step for GD and NGD. %strategies as a function of the gradient descent algorithm (CG versus NGD) and the choice of the loss function (KL versus W2). 
The simulation parameter values are set to $x_0^\star=1$, $\boldsymbol \varphi^* = \{ 0.5, 0.1,0.05 \}$, $\boldsymbol \varphi^{(0)} = \{ 1.5, 0.4,0.14\}$, $\sigma_\varepsilon = 0.05$, $N_\text{meas} = 41$, $t_{N_\text{meas}} = 2.2$, $\epsilon_{\text{KL}} = 10^{-2}$, $N_T = 13550$, $N_B= 29282$, $N_I  = 6655$, $N_R = 248897$.}
\label{fig:metapars_case3}
\end{figure}

\begin{figure}[htbp]
\begin{subfigure}{.5\textwidth}
 \centering
  % include first image
  \includegraphics[width=\linewidth]{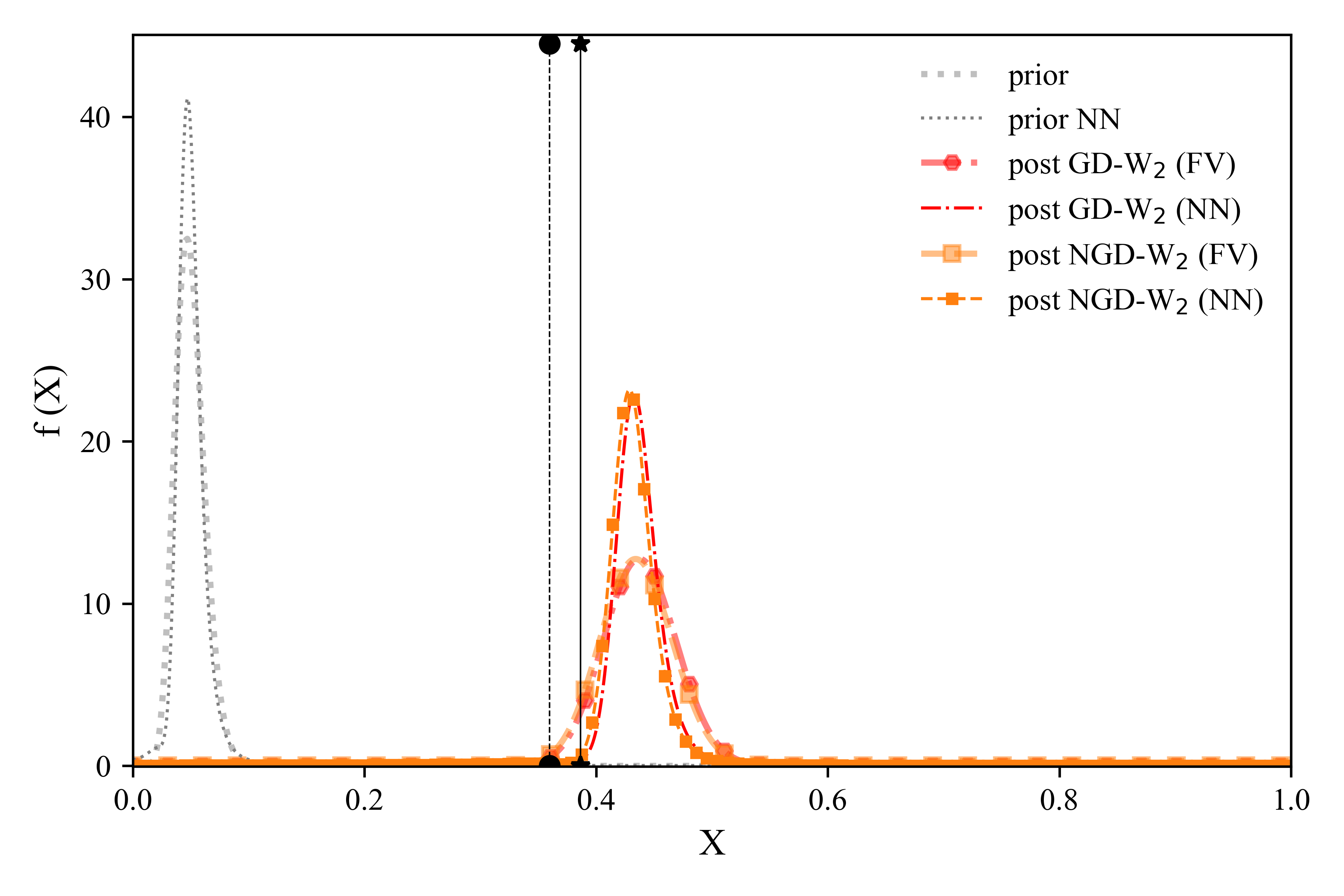}  
\end{subfigure}
\begin{subfigure}{.5\textwidth}
 \centering
  % include first image
  \includegraphics[width=\linewidth]{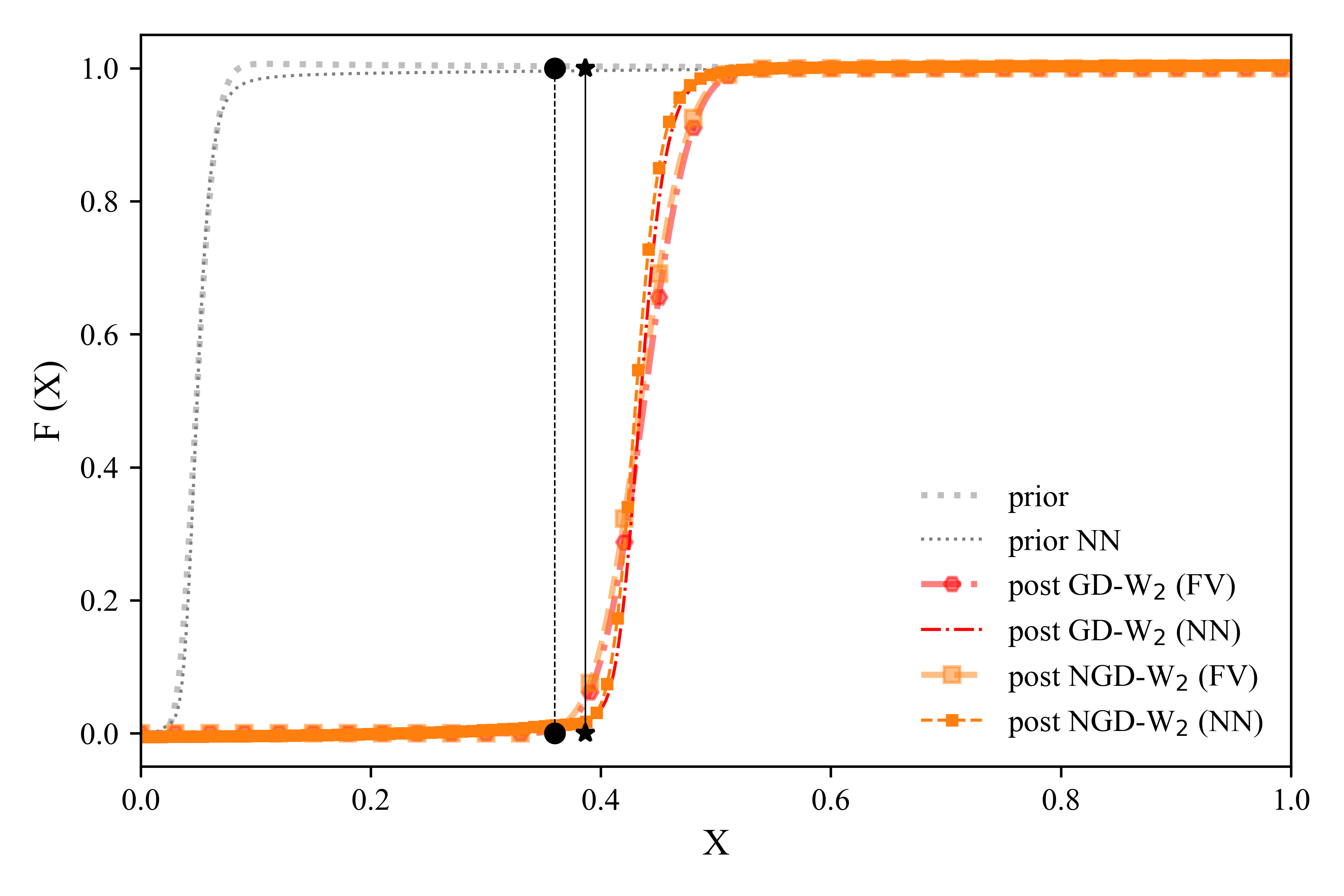}  
\end{subfigure}
\caption{Prior and posterior distributions at time $t_{N_{\text{meas}}}$ (both PDFs, on the left, and CDFs, on the right) obtained via GD-W$_2$ and NGD-W$_2$ minimization for the Langevin equation with colored noise. For each distribution, both the FV solution and NN approximation are shown. Black stars and circles mark the exact value $x(t_{N_\text{meas}})$ and its noisy observation $x_{N_\text{meas}}$, respectively. The simulation parameter values are the same as in \cref{fig:metapars_case3}.
}
\label{fig:posteriors_case3}
\end{figure}

In \cref{fig:posteriors_case3}, we present the posterior PDF and CDF of the state $x(t)$ at the final assimilation time $t_{N_\text{meas}}$. The CDF is computed as a FV solution of the CDF equation with   meta-parameters $\boldsymbol\varphi^{(N_\text{meas})}$, and the PDF as its derivative. Observations $\hat{\mathbf x}$ are assimilated, alternatively, via the GD-W$_2$ and NGD-W$_2$ optimization strategies. Both approaches yield posterior distributions that are close to the true state, with negligible differences between NGD-W$_2$ and GD-W$_2$. The use of the FV solution of the CDF equation leads to a slightly wider posterior than the reliance on its DNN surrogate does, possibly because of numerical diffusion. %The latter approach yields a posterior distribution that is  significantly closer to the true state; the use of the FV solution of the CDF equation leads to the wider posterior than the reliance on its DNN surrogate does. \textbf{[Yes?]} The accuracy of the GD posterior is compromised by convergence fluctuations at the end of the assimilation period (see \cref{fig:metapars_case3}), which also widen the posterior. 

Although not shown here, we found the KL- and $W_2$-based loss functions at %(\cref{fig:lossfunctionsKL_case3} and \cref{fig:lossfunctionsW2_case3}, respectively); cross-sections are shown for the loss function corresponding to 
at the first and later assimilation steps %($m=1$). The topology is 
to be smooth and not significantly different from each other. Yet, similar to the example in \cref{sec:problem2}, the differences are sufficient to prevent convergence in the KL case for poor choices of the prior. 

\section{Discussion and Conclusions}
\label{sec:concl}

We presented an information-geometric implementation of DA-MD, which yields computationally efficient data assimilation and parameter estimation for nonlinear problems with non-Gaussian system states. The forecast step is performed by employing the MD, an uncertainty propagation technique that yields a deterministic evolution equation for the CDF (or, equivalently, the PDF) of the state. This equation maps a set of meta-parameters (statistical properties of the random inputs) onto the system-state distribution, and defines a parameter space for a dynamic manifold of distributions. The analysis step is performed on this statistical manifold; it is formulated as sequential minimization of the discrepancy between an observational distribution and a predictive posterior distribution obeying the CDF equation with unknown (posterior) parameters. The observational PDF is the Bayesian posterior obtained as the product of the data model (i.e., the likelihood function) and the prior distribution obeying the CDF equation with the parameters from the previous assimilation step. 

Reliance on statistical discrepancy measures---the Kullback-Leibler divergence and the $L_2$ Wasserstein distance---confers exploitable geometric properties to the manifold of distributions. Specifically, it enables the use of NGD, an efficient optimization technique. Our numerical experiments revealed the $W_2$-based DA-MD to be more robust to the choice of a prior than its KL-based counterpart. 

For one-dimensional (univariate) distributions, $W_2$ is defined in terms of system-state CDFs, and KL in terms of corresponding PDFs. This argues in favor of the $W_2$-based DA-MD, since CDFs are smoother and numerical solution of CDF equations is easier. This facilitates the use of invertible DNNs as a surrogate model in the probabilistic space to facilitate and accelerate optimization and calculation of the geometric metric tensors.

Future work will focus on the identification of ambiguity sets and their dynamics on statistical manifolds \cite{boso-2020-ambiguity}, their evolution and their update with observations. We also plan to explore the use of different data models, the impact of alternative parameterizations of a statistical manifold on DA-MD performance, and the latter's implications for sensitivity analysis.

%\ethics{Insert ethics text here.}

% \dataccess{There are no data sharing issues since all of the numerical information is provided in
% the figures produced by solving the equations in the paper.}

% \aucontribute{FB conceived the study, developed the analysis, carried out the simulations, and
% drafted the manuscript; DT participated in the design of the study, participated in the analysis, critically
% revised the manuscript. All authors gave final approval for publication and agree to be held accountable for the work performed therein.}

% \competing{The authors declare no competing interests.}

% This work was supported in part by Air Force Office of Scientific Research under award number
% FA9550-18-1-0474, and by a gift from TOTAL.

%\ack{Insert acknowledgment text here.}

%\disclaimer{Insert disclaimer text here.}

\appendix

\section{CDF equation for the stochastic ODE} 
\label{app:MD}

We summarize the MD for the three test problems from \cref{sec:problem}. The original derivations can be found in \cite{tartakovsky2016method}, \cite{jazwinski1970stochastic} and~\cite{maltba-2018-nonlocal}, respectively. The first two results are exact, whereas the third one is approximate and has been verified against Monte Carlo simulations in \cite{barajas2016probabilistic,maltba-2018-nonlocal}.

\subsection{Stochastic ODE with random initial conditions}
\label{sec:case1}

We consider \eqref{eq:ODEphysical} with a smooth deterministic function $s(x,t;\boldsymbol \theta)$; random initial state $x_0 \in \mathbb R$ is described by a given CDF $F_0(X;\boldsymbol \varphi_0)$ with statistical parameters $\boldsymbol \varphi_0$. To derive an equation for $F(X;t)$, the CDF of $x(t)$,  we first define the raw CDF $\Pi = \mathcal H(X - x(t))$ whose ensemble mean is $\langle \Pi(X;t) \rangle = F(X;t)$. Next, we multiply \eqref{eq:ODEphysical} by $- \partial \Pi / \partial X$ and use the properties of the Heaviside function $\mathcal H(\cdot)$ to obtain
\begin{align}\label{eq:Pi1}
    \frac{\partial \Pi}{\partial t} + s(X,t; \boldsymbol \theta ) \frac{\partial \Pi}{\partial X} = 0.
\end{align}
Since $s(\cdot)$ is deterministic, the ensemble average of~\eqref{eq:Pi1} yields \eqref{eq:CDF1}, 
which is a special case of~\eqref{eq:CDFeq_closed} with $\mathcal U = s(X,t; \boldsymbol \theta)$ and $\mathcal D = 0$. %$\Pi(X,t=0) = \Pi_0(X) = \mathcal H(X - x_0)$.

The derivation of a corresponding PDF equation starts with the definition of a raw PDF $\pi = \delta(X-x(t))$, where $\delta(\cdot)$ is the Dirac distribution. A procedure similar to above yields  \cite[sec. 2.1]{tartakovsky2016method} 
\begin{align} \label{eq:PDF1}
 \frac{\partial f}{\partial t} + \frac{\partial s(X,t;{\boldsymbol \theta}) f}{\partial X} = 0. 
\end{align}
This equation can also be obtained by differentiation of \eqref{eq:CDF1} with respect to $X$. %$f(X,t=0) = f_0(X;\boldsymbol \varphi_0) = \frac{\partial F_0(X;\boldsymbol \varphi_0)}{\partial X}$

\subsection{MD for the Langevin equation with white noise}

Consider a Langevin equation, \eqref{eq:ODEphysical} with $s(x; w) \equiv s_d(x,t) + s_w(x,t) w(t)$ where $w(t)$ is a white standard Gaussian process (with zero mean and unit variance). The deterministic functions $s_d$ and $s_w$ are such that $s(x;w)$ is integrable with respect to $t$ in the mean square sense \cite[Sec. 4.1]{jazwinski1970stochastic}.  The derivation of a PDF equation for $x(t)$ is relatively straightforward, and leads to the Fokker-Planck equation (a.k.a. Kolmogorov's forward equation)~\cite[Sec. 4.9]{jazwinski1970stochastic}
\begin{equation}\label{eq:FPE}
     \frac{\partial f}{\partial t} + \frac{\partial \, s_d(X,t) f}{\partial X} = \frac{1}{2} \frac{\partial^2 s_w^2(X,t) f}{\partial X^2},
\end{equation}
It is formally valid if $f(X;t)$ is well-behaved at infinity, and is subject to initial and boundary conditions condition $f(X; 0) = f_0(x)$ and  $f( \pm \infty; t) = 0$.

An equivalent CDF version of the Fokker-Planck equation \eqref{eq:FPE} can be obtained via integration of \eqref{eq:FPE} over $X \in \Omega$ 
\begin{equation}\label{eq:FPE_CDF}
    \frac{\partial F}{\partial t} + s_d(X,t) \frac{\partial F}{\partial X}  = \frac{1}{2} \frac{\partial }{\partial X} \left( s_w^2(X,t) \frac{ F}{\partial X} \right),
\end{equation}
subject to $F(x;0) = F_0(X)$, $F(X_{\text{min}},t) = 0$, and $F(X_{\text{max}},t) = 1$.

In~\eqref{eq:ODE2}, $s(x; w) = - a(t) x(t)$ where the random process $a(t)$ has the constant mean $\mu_a$ and standard deviation $\sigma_a$. This translates into $s_d(x,t) = - \mu_a x$ and $s_w = - \sigma_a x$, so that the  coefficients $\mathcal U$ and $\mathcal D$ in~\eqref{eq:CDFeq_closed} become $\mathcal U = - \mu_a X$ and $\mathcal D = (\sigma_a^2 /2) X^2$, with $\boldsymbol \varphi = \{\mu_a,\sigma_a \}$. 

\subsection{MD for the Langevin equation with colored noise}
\label{app:problem3}

Consider~\eqref{eq:ODEphysical} with $s(x; w) \equiv -a(t) x(t)$, where $a(t) = \mu_a + w(t)$ and $w(t)$ is a correlated standard Gaussian process.  The MD for stochastic/random (Langevin) ODEs with temporally correlated forcings requires closure approximations. These include the semi-local approximation~ \cite{barajas2016probabilistic,maltba-2018-nonlocal}, which compares favorably with Monte Carlo simulations and a local closure approximation in terms of both accuracy and computational efficiency. For the sake of completeness, we summarize the derivation of the PDF equation and its semi-local closure approximation for the specific form of the Langevin equation described above. We start by deriving an equation for the raw PDF $\pi(X,t) = \delta (X - x(t))$, whose ensemble mean is the PDF, $f(X;t) = \langle \pi \rangle$. Multiplying our ODE by $- \partial \pi / \partial X$ and using the properties of the Dirac delta function $\delta(\cdot)$, we obtain
\begin{align}\label{eq:app-pi}
   \frac{\partial \pi }{ \partial t } + a(t) \frac{\partial \pi }{ \partial X} = 0. \qquad
  s(X,t) = \langle s(X,t) \rangle + s'(X,t;w); \quad \langle s \rangle = - \mu_a X, s'= - w(t) X
\end{align} 
We use the Reynolds decomposition $\mathcal A = \langle \mathcal A \rangle + \mathcal A'$ to represent relevant random processes $\mathcal A$ as the sums of their ensemble means $\langle \mathcal A \rangle$ and zero-mean fluctuations around these means, $\mathcal A'$. Since $\pi = f + \pi'$, taking the ensemble mean of this equation yields an unclosed equation for the PDF $f(X;t)$,
\begin{equation}\label{eq:app-unclosed}
    \frac{\partial f}{\partial t} + \mu_a \frac{ \partial f}{\partial X} + \frac{\partial \langle w'(t) \pi'(X,t) \rangle}{\partial X} = 0, \qquad \text{subject to} \quad f(X;0) = f_0.
\end{equation}
A closure approximation is needed to render the cross-correlation term $\langle w'(t) \pi'(X,t) \rangle$ computable. Subtracting~\eqref{eq:app-unclosed} from~\eqref{eq:app-pi}, we obtain an equation for random fluctuations $\pi'(X,t)$, 
\begin{equation}\label{eq:piprimecolored}
    \frac{\partial \pi'}{\partial t} + \mu_a \frac{ \partial \pi'}{\partial X} = \frac{\partial ( \langle s'(X,t) \pi'(X,t) \rangle - s' \pi ) }{\partial X}, \qquad \text{subject to} \quad \pi'(X,t=0) = 0.
\end{equation}
The deterministic Green's function for~\eqref{eq:piprimecolored}, $G(X,t;\Xi,\tau)$, is a solution of 
\begin{align}\label{eq:greencolored}
    \frac{\partial G}{\partial \tau} + \mu_a \frac{\partial G}{\partial \Xi} = - \delta (X-\Xi) \delta(t-\tau)
\end{align}
with homogeneous initial (at $\tau=t$) and boundary conditions at infinity. Its analytical solution, obtained, e.g., via the method of characteristics, is $G(X,t;\Xi,\tau) = \mathcal H(t-\tau) \delta(X - \Xi \exp(- \mu_a (t-\tau)))$. Hence, the path-wise solution of \eqref{eq:piprimecolored} is
\begin{equation}\label{eq:piprime}
    \pi'(X,t) = \int_0^t \int_{-\infty}^\infty G(X,t;\Xi,\tau) \frac{\partial}{\partial \Xi} \left[ \langle w' (\Xi,\tau) \pi' (\Xi,\tau)\rangle - w'(\Xi,\tau) \pi(\Xi,\tau) \right] \text d \tau \text d \Xi.
\end{equation}
A closure approximation for $\langle w'(t) \pi'(X,t) \rangle$ is constructed by multiplying~\eqref{eq:piprime} with $w'(t)$, taking the ensemble mean, and neglecting the third-order correlation term,
\begin{align}\label{eq:app-last}
    \langle w'(X,t) \pi'(X,t) \rangle = - \int_0^t \int_{-\infty}^\infty G(X,t;\Xi,\tau) \frac{\partial}{\partial \Xi} \left( C_w (X,t;\Xi,\tau) f(\Xi,\tau)\right) \text d \Xi \text d \tau,
\end{align}
where $C_w(X,t;\Xi,\tau) = \langle w'(X,t) w'(\Xi,\tau) \rangle$ is the auto-covariance of the random noise $w(t)$. Substituting this expression into~\eqref{eq:app-unclosed} yields a nonlocal (integro-differential) PDF equation. Accounting for the analytical expression for $G$,~\eqref{eq:app-last} is approximated semi-locally as
\begin{align}
    \langle w'(X,t) \pi'(X,t) \rangle = - X f(X,t) \int_0^t C_w (t,\tau) \text d \tau  - X^2 \frac{\partial f(X,t)}{\partial X} \int_0^t C_w (t,\tau) \text d \tau.
\end{align}
This yields the closed CDF equation \eqref{eq:CDFeq_closed} with \eqref{eq:UDtest3}. If $w(t)$ were white noise, i.e., if $C_w (t,\tau) = \delta(t-\tau)$, then the resulting PDF equation would reduce to the Fokker-Planck equation.

\end{document}